\documentclass{amsart}
\usepackage{amsfonts, amsbsy, amsmath, amssymb}
\newtheorem{thm}{Theorem}[section]
\newtheorem{lem}[thm]{Lemma}
\newtheorem{cor}[thm]{Corollary}

\newtheorem{exmp}[thm]{Example}

\newtheorem{conj}[thm]{Conjecture}
\numberwithin{equation}{section}

\theoremstyle{definition}
\newtheorem{defn}{Definition}[section]
\begin{document}

\title[Affinity of a Permutation]{The Affinity of a Permutation of \\
a Finite Vector
Space}

\author{W. Edwin Clark}
\address{Department of Mathematics,
University of South Florida, Tampa, FL 33620}
\email{eclark@math.usf.edu}

\author[Xiang-dong Hou]{Xiang-dong Hou*}
\address{Department of Mathematics,
University of South Florida, Tampa, FL 33620}
\email{xhou@tarski.math.usf.edu}

\author{Alec Mihailovs}
\address{Department of Mathematics,
Tennessee Technological University, Cookeville, TN 38505
}
\email{alec@mihailovs.com}

\thanks{* Partially supported by NSA grant MDA 904-02-1-0080}


 
\begin{abstract} 
For a permutation $f$ of an $n$-dimensional vector
space $V$ over a finite field of order $q$ we let $k\text{-affinity}(f)$
denote the number of $k$-flats $X$ of $V$ such that $f(X)$ is also a
$k$-flat. By $k\text{-spectrum}(n,q)$ we mean the set of integers 
$k\text{-affinity}(f)$ where
$f$ runs through all permutations of $V$. The problem of the complete
determination of $k\text{-spectrum}(n,q)$ seems very difficult except for small or special values of the parameters. However, we are able to establish that
$0 \in k\text{-spectrum}(n,q)$ in the following cases: (i) $q\ge 3$ and 
$1\le k\le n-1$; (ii) $q=2$, $3\le k\le n-1$; (iii) $q=2$, $k=2$, $n\ge 3$ odd. 
The maximum of $k\text{-affinity}(f)$ is,
of course, obtained when $f$ is any semi-affine mapping. We conjecture
that the next to largest value of $k\text{-affinity}(f)$ is when $f$ is a
transposition and we are able to prove this when $q=2$, $k=2$,  $n \ge 3$ and
when $q\ge 3$, $k=1$, $n \ge 2$. 
\end{abstract}
\medskip

\maketitle

%%%%%%%%% section 1 &&&&&&&&&&&&&&&&&&&&&&&&&&&&&&&&&&&&
\section{Introduction}

It is a classical result, see, \emph{e.g.}, Snapper and Troyer
\cite{Sna71}, that if
$V$ is an $n$-dimensional vector space over a
field $F$ such that $n\ge 2$ and $|F| \ge 3$ then a bijection $f: V \to V$ which
takes 1-flats to 1-flats is a semi-affine  mapping, that is, there is
an automorphism $\sigma$ of $F$, an additive automorphism $g: V \to V$
and a vector $b \in V$ such that $g(\alpha x) = \sigma(\alpha)g(x)$
for all $x \in V$,  $\alpha  \in F$ and 
\[
f(x) = g(x) + b \qquad \mbox{for all $x \in V$}.
\]
We remark that if the automorphism $\sigma$ is the identity then $g$ is just a
non-singular linear mapping and $f$ is said to be affine. This will be the case
when $F$ has no non-trivial automorphisms.

The above result is not true when $|F| = 2$. In this case, a 1-flat in $V$ is just a two element subset, 
hence every permutation of $V$
takes all 1-flats to 1-flats.
However, the above result has an easy analog for the case $|F| = 2$:
A permutation of $V$ which takes every 2-flat to a 2-flat must be
affine (cf. \cite{Hou}).   

Let ${\mathbb F}_q$ be the finite field with $q$ elements and let ${\mathbb F}_q^n$ be the $n$-dimensional vector space over ${\mathbb F}_q$.
In this paper, we are
concerned with permutations of ${\mathbb F}_q^n$. Let $\text{Per}({\mathbb F}_q^n)$ denote the group of all permutations of ${\mathbb F}_q^n$.
Recall that a $k$-flat (or $k$-dimensional
affine subspace) $X$ in ${\mathbb F}_q^n$ is a coset $U + x$ of a $k$-dimensional subspace $U$ of ${\mathbb F}_q^n$.

\begin{defn} 
For $f \in \text{Per}({\mathbb F}_q^n)$ and $0 \le k \le n$ we
define $k\text{-affinity}(f)$ to be the number of $k$-flats $X$ in ${\mathbb F}_q^n$ such that 
$f(X)$ is a $k$-flat. 
We define $k\text{-coaffinity}(f)$ to be the number of $k$-flats $X$ in 
${\mathbb F}_q^n$
such that $f(X)$ is not a $k$-flat.
\end{defn}
 
It is well known that the number of $k$-dimensional subspaces of
${\mathbb F}_q^n$ is given by the $q$-binomial coefficient
\[
{\genfrac{[}{]}{0pt}{0}n k}_q =
\frac{(q^n-1)(q^{n-1}-1)\dots(q^{n-k+1}-1)}{(q^k-1)(q^{k-1}-1)\dots(q^1-1)}
\]
and the number of $k$-flats in ${\mathbb F}_q^n$ is given by
\[
q^{n-k}{\genfrac{[}{]}{0pt}{0}n k}_q.
\]
It follows that 
\[
k\text{-affinity}(f) + k\text{-coaffinity}(f) = q^{n-k}{\genfrac{[}{]}{0pt}{0}n k}_q
\]
for all permutations $f$ of ${\mathbb F}_q^n$ and all $0 \le k \le n$.

The cases $k = 0$ and $k = n$ are trivial and we shall ignore them.

\begin{defn} For integers $0 \le k \le n$ and prime power $q$, we
define $k\text{-spectrum}(n,q)$ to be the set of values 
$k\text{-affinity}(f)$ for all $f \in \text{Per}({\mathbb F}_q^n)$.
\end{defn} 

The present paper is a continuation of the second author's work \cite{Hou}. In \cite{Hou}, the notion of 2-affinity of permutations of 
${\mathbb F}_2^n$ was implicitly introduced and permutations of 
${\mathbb F}_2^n$ with 2-affinity 0 were studied. We point out that a
permutation $f\in \text{Per}({\mathbb F}_2^n)$ with 
$2\text{-affinity}(f)=0$ is an \emph{almost perfect
nonlinear} (APN) permutation. APN permutations arose in cryptography as
a means to resist the differential cryptanalysis \cite{Car98, Nyb95}. APN permutations of ${\mathbb F}_2^n$ are known to exist for odd $n\ge 3$ (\cite{Car98, Nyb94}) and not to exist for $n=4$ (\cite{Hou}). Their 
existence for even $n\ge 6$ is an open question. For recent work on APN permutations and related topics, we refer the reader to \cite{Bet94, Car98, Dob99, Dob01, Hou}. However, we must remind the reader that this paper
is not a response to any problem from cryptography. Rather, it is a pure mathematical exploration.

Our primary interest is the set $k\text{-spectrum}(n,q)$. In particular,
we would like to know if $0\in k\text{-spectrum}(n,q)$ and what the second largest number in $k\text{-spectrum}(n,q)$ is. (The largest number in
$k\text{-spectrum}(n,q)$ is, of course, $q^{n-k}{\genfrac{[}{]}{0pt}{1}n k}_q$.) In Section 2, we show that with few exceptions, $0\in 
k\text{-spectrum}(n,q)$. The result of Section 2 relies on an inequality involving $q$-binomial coefficients whose proof is given in Section 3.
Hou \cite{Hou} showed that 
2-spectrum(4,2) is 
\[
\{\text{5 -- 20, 22, 24 -- 26, 28, 30, 32, 36, 38, 44, 48, 52, 56, 76, 84, 140} \}
\]
where $a$ -- $b$ denotes all integers from $a$ to $b$. 
More examples of $k$-spectra are given in Section 4. 
In Section 5, we determine $(n-1)\text{-spectrum}(n,2)$ completely.
These examples and results led to the
conjecture that the next to largest $k$-affinity is that of a transposition. We
compute the $k$-affinity $T(n,k,q)$ of a transposition in 
$\text{Per}({\mathbb F}_q^n)$ in Section 6. We call this conjecture \emph{The Threshold Conjecture}  since it says that if
$k\text{-affinity}(f)>T(n,k,q)$  then $f$ takes every $k$-flat to a $k$-flat.
We prove that the conjecture holds for $q=2$, $k=2$,
$n\ge 3$ in Section 7 and for $q > 2$, $k=1$, $n \ge 2$ in Section 8.
 
%%%%%%%%%%%%%%%%% section 2 %%%%%%%%%%%%%%%%%%%%%%%%%%
\section{When $k\text{-affinity}(f)=0$}

It should be noted that there appears to be no clear relationship between
$k$-affinity($f$) and $\ell$-affinity($f$). For example, there are
permutations $f_1,f_2,f_3,f_4$ in $\text{Per}({\mathbb F}_3^3)$ such that 
\[
\begin{split}
& 1\text{-affinity}(f_1)=1\  \text{and}\  2\text{-affinity}(f_1)=0\cr 
& 1\text{-affinity}(f_2)=0\  \text{and}\  2\text{-affinity}(f_2)=1\cr 
& 1\text{-affinity}(f_3)=0\  \text{and}\  2\text{-affinity}(f_3)=0\cr 
& 1\text{-affinity}(f_4)=1\  \text{and}\  2\text{-affinity}(f_4)=1\cr 
\end{split}
\]

In the following theorem, we see that with few exceptions there is a permutation $f\in \text{Per}({\mathbb F}_q^n)$ such that simultaneously 
$k\text{-affinity}(f)=0$ for all $1\le k\le n-1$.

\begin{thm} \label{T2.1} 
\
\begin{enumerate}
\item[(i)] If $q = 2$ and $n\ge 3$ is odd, there exists $f \in
\text{\rm Per}({\mathbb F}_2^n)$ such that $2\text{\rm-affinity}(f)= 0$.

\item[(ii)] If $q = 2$ and $n \ge 4$ , there exists $f \in \text{\rm Per}({\mathbb F}_2^n)$ such that $k\text{\rm-affinity}(f)= 0$ for all $3 \le k \le n-1$.

\item[(iii)] If $q \ge 4$ and $n \ge 2$, there exists $f \in
\text{\rm Per}({\mathbb F}_q^n)$ such that $k\text{\rm-affinity}(f)= 0$ for all $1 \le k \le n-1$.

\item[(iv)] If $q = 3$ and $n \ge 3$, there exists $f \in
\text{\rm Per}({\mathbb F}_3^n)$ such that $k\text{\rm-affinity}(f)= 0$ for all $2 \le k \le n-1$.

\item[(v)] If $q = 3$ and $n\ge 2$, there exists $f \in
\text{\rm Per}({\mathbb F}_3^n)$ such that $1\text{\rm-affinity}(f)= 0$.
\end{enumerate} 
\end{thm}

The proof of Theorem~\ref{T2.1} is spread out in parts in the rest of this section. Part (i) of Theorem~\ref{T2.1} is well known. (See \cite{Dob99, Dob01}
for several families of permutations of ${\mathbb F}_2^{2m+1}$ with
2-affinity 0.) It also follows from the following example in which 
we compute the 2-affinity of the permutation $f$ of ${\mathbb F}_2^n$
($\cong{\mathbb F}_{2^n}$) defined by $f(x)=x^{2^n-2}$. We remark that this permutation has been discussed by Nyberg \cite{Nyb94} in terms of differential uniformity and that our computation is slightly different
from that of \cite{Nyb94}.

\begin{exmp}\label{E2.2}
\rm
Identify ${\mathbb F}_2^n$ with ${\mathbb F}_{2^n}$. Define $f\in \text{\rm Per}({\mathbb F}_{2^n})$ by $f(x) = x^{2^n-2}$ for
$x\in {\mathbb F}_{2^n}$. Note that 
\[
f(x) =
\begin{cases}
0 &\text{if}\ x = 0, \cr
\frac 1x &\text{if}\  x \ne 0. \cr
\end{cases}      
\]
We claim that
\[
2\text{\rm-affinity}(f)=
\begin{cases}
0 &\text{if $n$ is odd},\cr
\frac{2^n-1}{3}&\text{if $n$ is even}.\cr
\end{cases}  
\]

Suppose $X\subset{\mathbb F}_{2^n}$ is some  2-flat such that $f(X)$ is also a 2-flat. Then we can write
$X =\{x,y,z,w\}$ where $x+y+z+w = 0$. Suppose first that $x$, $y$, $z$, and $w$ are
all nonzero. Then we have
$f(X) =\{\frac 1x,\frac1y,\frac1z,\frac1w\}$ and
$\frac1x+\frac1y+\frac1z+\frac1w = 0$. It follows that $w = x+y+z$ and  
\[
0=\frac 1 x + \frac 1 y +\frac 1 z + \frac 1 {x+y+z} = \frac {(x+y)(x+z)(y+z)}{xyz(x+y+z)}.
\]
It follows that $x = y$ or $x = z$ or $y = z$ which contradicts the assumption
that $X$ is a 2-flat. Thus without loss of generality we may assume that $w = 0$.
Then $X = \{0,x,y,x+y \}$,  $f(X) = \{0, \frac1x,\frac1y,\frac1{x+y}\}$ and 
\[
\frac 1 x + \frac 1 y +\frac 1 {x+y} = 0.
\]
This is equivalent to 
\begin{equation}\label{2.1}
\left(\frac y x \right)^2 + \left(\frac y x \right) + 1  = 0.
\end{equation}
Hence $\frac yx$ is a root of the irreducible polynomial $g(t)=t^2+t+1 \in \mathbb{F}_2[t]$. Therefore
${\mathbb F}_2(\frac yx)$ is a subfield of $\mathbb{F}_{2^n}$ 
with $[{\mathbb F}_2(\frac yx):{\mathbb F}_2]=2$. It follows
that $n$ is even. So if
$n$ is odd, no such $x$ and $y$ exist and $2\text{\rm-affinity}(f) = 0$.

On the other hand, if $n$ is even, $g(t)$ has two roots $\beta$ and $1+\beta$ in
${\mathbb F}_{2^n}$ and $K=\{0,1,\beta,1+\beta\}$ is the unique subfield of order 4
in ${\mathbb F}_{2^n}$. It follows from (\ref{2.1}) that $\frac yx =\beta$ or $1+\beta$.
In both cases,
\begin{equation*}
\begin{split}
   X &=\{ 0,x, y , x + y \} \cr
    &=\{ 0,x,x\beta , x(1 + \beta ) \} \cr
    &=x K \cr
    &=\{ 0 \} \cup x K^* \cr
\end{split}
\end{equation*}
where $K^*$ is the multiplicative group of $K$ and $x K^*$ is a coset of the subgroup
$K^*$ in ${\mathbb F}_{2^n}^*$. There are $(2^n-1)/3$ cosets of $K^*$ in
${\mathbb F}_{2^n}^*$ and therefore the same number of 2-flats of the form $xK$.
Since $f(x K) = \frac 1 xK$, it follows that $2\text{\rm-affinity}(f)=
(2^n-1)/3$, as claimed.

Note that when $n=4$, $\frac{2^4-1}3 =5$
is the  minimum 2-affinity of permutations of ${\mathbb F}_2^4$ (\cite{Hou}).
\end{exmp}

The proof of parts (ii) -- (iv) of Theorem~\ref{T2.1} relies on the
following theorem whose proof will be given in Section 3.

\begin{thm}\label{T2.3} 
Let $q,\ m,\ n$ be integers such that $n>m$ and $q,\ m$ satisfy one of the following conditions:
\begin{enumerate}
\item[(a)] $q=2$, $m=3$,

\item[(b)] $q=3$, $m=2$,

\item[(c)] $q\ge 4$, $m=1$.
\end{enumerate}
Then
\begin{equation}\label{2.2}
\sum_{k=m}^{n-1}q^{2(n-k)}{\genfrac{[}{]}{0pt}{0}n k}_q^2 q^k!\; (q^n-q^k)! < q^n! .
\end{equation}
\end{thm}

\begin{proof}[Proof of Theorem~\ref{T2.1} (ii) -- (iv)]
Let $\Phi_k$ denote the set of all $f \in
\text{\rm Per}({\mathbb F}_q^n)$ such that $k\text{-affinity}(f)\ge 1$. Let $X$ be any fixed $k$-flat
in ${\mathbb F}_q^n$ and define 
\[
S_X = \{ f \in \text{\rm Per}({\mathbb F}_q^n) :  f(X) = X \}.
\]
Note that since $|X| = q^k$ we have $|S_{X}| = q^k!\,(q^n-q^k)!$. 

The group of invertible affine transformations of ${\mathbb F}_q^n$,
i.e., the general affine group ${\rm AGL}(n,{\mathbb F}_q)$, acts transitively on the set of all $k$-flats of ${\mathbb F}_q^n$. Hence for
every $k$-flat $W$, there exists $\alpha_W\in {\rm AGL}(n,{\mathbb F}_q)$
such that $\alpha_W(W)=X$. Let $f \in \Phi_k$ and assume that $f(W) = Z$ where $W$ and
$Z$ are $k$-flats. Then $\alpha_Z\circ f\circ\alpha_W^{-1}\in S_X$, i.e.,
$f\in\alpha_Z^{-1}\circ S_X\circ\alpha_W$.
Since there are $q^{n-k}{\genfrac{[}{]}{0pt}{1}n k}_q$
$k$-flats in ${\mathbb F}_q^n$ it follows that 
\[
|\Phi_k| \le  \left(q^{n-k} \genfrac{[}{]}{0pt}{0}n k_q \right)^2
|S_X|=q^{2(n-k)}\genfrac{[}{]}{0pt}{0}n k_q^2  q^k!\,(q^n-q^k)!.
\]
Hence if $q$ and $m$ satisfy one of the conditions in Theorem~\ref{T2.3}
and $n>m$, we have
\[
\sum_{k=m}^{n-1}|\Phi_k|\le \sum_{k=m}^{n-1}
q^{2(n-k)}\genfrac{[}{]}{0pt}{0}n k_q^2  q^k!\,(q^n-q^k)! < q^n!.
\]
Thus there exists $f\in \text{\rm Per}({\mathbb F}_q^n)$ such that $f\notin
\bigcup_{k=m}^{n-1}\Phi_k$.
\end{proof}

Note that inequality (\ref{2.2})
does not cover the case $q = 3$ and $m = 1$, that is, part (v) of 
Theorem~\ref{T2.1}. This case is dealt with as a corollary to the following lemma.

\begin{lem} \label{L2.4} 
Let $F$ be any field. If 
the permutations  $f: F^n \to F^n$ and $g:F^m \to F^m$ each have
1-affinity 0 then the permutation $f \times g:F^n \times F^m \to F^n \times F^m$ has 
1-affinity 0.
\end{lem}

\begin{proof} 
Assume to the contrary that there exists a 1-flat $X$ in $F^n\times F^m$
such that $(f\times g)(X)$ is also a flat. Let $\pi_1:F^n\times F^m\rightarrow F^n$  and $\pi_2:F^n\times F^m\rightarrow F^m$ be the projections. Then either $\pi_1(X)$ is a 1-flat in $F^n$ or $\pi_2(X)$
is a 1-flat in $F^m$. Without loss of generality, assume that the former
is the case. Thus $f\bigl(\pi_1(X)\bigr)=\pi_1\bigl((f\times g)(X)\bigr)$
is a 1-flat in $F^n$, which is impossible since $1\text{-affinity}(f)=0$.
\end{proof}

\begin{cor}\label{C2.5}
If $n\ge 2$, there exists $f \in \text{\rm Per}({\mathbb F}_3^n)$ such that
$1\text{\rm -affinity}(f)=0$.
\end{cor}

\begin{proof}
By Lemma~\ref{L2.4}, 
it suffices to show that ${\mathbb F}_3^2$ and
${\mathbb F}_3^3$ have permutations of 1-affinity 0. Such permutations are easily
 found through a computer search. For ${\mathbb F}_3^2$, label the elements
$(0,0),\, (0,1),\, \dots,\, (2,2)$ with $0,1,\dots,8$. A desirable permutation $f$ is given by 
\[
\bigl(f(0),\cdots,f(8)\bigr)= (0,\ 1,\ 8,\ 2,\ 3,\ 4,\ 5,\ 6,\ 7).
\]

For ${\mathbb F}_3^3$, we label the elements 
$(0,0,0),\, (0,0,1),\,\dots,\, (2,2,2)$ with  $0,1, \dots, 26$.
A desirable permutation $f$ in this case is given by
\[
\begin{split}
\bigl(f(0),\cdots,f(26)\bigr)=\,& (0,\ 1,\ 24,\ 2,\ 3,\ 4,\ 5,\ 6,\ 7,\ 8,\ 9,\ 10,\ 11,\ 12,\ 13,\cr 
&14,\ 25,\ 15,\ 16,\ 17,\ 26,\ 18,\ 19,\ 23,\ 20,\ 21,\ 22).\cr
\end{split}
\]
\end{proof}

%%%%%%%%%%%%%%%% section 3 %%%%%%%%%%%%%%%%%%%%%%%%%%%%%%%%%%%%%
\section{Inequalities between Binomial and $q$-Binomial Coefficients}

In this section we assume that $i,k,m,n,q$ are integers.

\begin{lem}\label{L3.1}
For $q>2$, $n>k\geq 1$ and $q=2$, $n>k\geq 2$, 
\begin{equation}\label{3.1}
\frac{\genfrac{[}{]}{0pt}{0}n k_q^2}{\dbinom{q^n}{q^k}} <
\frac{1}{q^{q^k-2k}}\cdot
\frac{\genfrac{[}{]}{0pt}{0}{n-1} k_q^2}{\dbinom{q^{n-1}}{q^k}} .
\end{equation}
\end{lem}

\begin{proof}
Inequality (\ref{3.1}) is equivalent to
\begin{equation}\label{3.2}
\frac{(q^n-1)^2}{q^n(q^n-1)\dots(q^n-(q^k-1))} <
\frac{(q^n-q^k)^2}{q^n(q^n-q)\dots(q^n-(q^k-1)q)},
\end{equation}
which can be further rewritten as
\[
\frac 1{\prod_{i=2}^{q^k-2}(q^n-i)}\cdot \frac{q^n-1}{q^n-(q^k-1)}<
\frac 1{\prod_{\substack{1\le i\le q^k-2\cr i\ne q^{k-1}} } (q^n-qi)}
\cdot\frac{q^n-q^k}{q^n-(q^k-1)q}.
\]
Clearly,
\[
\prod_{i=2}^{q^k-2}(q^n-i)\ge \prod_{\substack{1\le i\le q^k-2\cr i\ne q^{k-1}} } (q^n-qi).
\]
Thus it suffices to show that
\begin{equation}\label{3.3}
\frac{q^n-1}{q^n-(q^k-1)}<\frac{q^n-q^k}{q^n-(q^k-1)q}.
\end{equation}
Inequality (\ref{3.3}) follows from  
\[
\begin{split}
&(q^n-q^k)(q^n-(q^k-1))-(q^n-1)(q^n-(q^k-1)q)\cr
=\,&(q^k-1)(q^n(q-2)+q^k-q)\cr
>\,&0.\cr
\end{split}
\]
\end{proof}

\begin{lem}\label{L3.2}
For $q\geq 4$, $k\geq 1$, or $q=3$, $k\geq 2$, or $q=2$, $k\geq 3$,
\begin{equation}\label{3.4}
\frac{\genfrac{[}{]}{0pt}{0}{k+1} k_q^2}{\dbinom{q^{k+1}}{q^k}} <
\frac{1}{q^{q^k-k}}.
\end{equation}
\end{lem}

\begin{proof}
The left hand side of (\ref{3.4}) equals 
\[
\frac{q^{k+1}-1}{q^{k+1}}\cdot 
\frac{2\cdot 3}{(q-1)^2(q^{k+1}-q^k+1)(q^{k+1}-q^k+2)}\cdot q^k\cdot 
\prod_{i=1}^{q^k-4}\frac{q^k-i}{q^{k+1}-i-1}.
\]
In this product, 
\[
\frac{q^{k+1}-1}{q^{k+1}}<1 
\]
and for every $1\le i\le q^k-4$,
\begin{equation}\label{3.5}
\frac{q^k-i}{q^{k+1}-i-1}\leq \frac{1}{q}
\end{equation} 
(To see (\ref{3.5}), note that since $q\ge 2$, we have
$q^{k+1}-i-1\geq q^{k+1}-qi$.)
Therefore, it suffices to show that 
\[
\frac{6}{(q-1)^2(q^{k+1}-q^k+1)(q^{k+1}-q^k+2)}\le \frac 1{q^4}.
\]
Let
\[
f(q,k)=(q-1)^2\Bigl(q^{k-2}(q-1)+\frac 1{q^2}\Bigr)\Bigl(q^{k-2}(q-1)+\frac 2{q^2}\Bigr).
\]
It suffices to show that
\[
f(q,k)\ge 6.
\]

The function
$f(q,k)$ is increasing with respect to $k$ for fixed $q>1$. For $k\ge 2$ and $q\ge 2$ or $k=1$ and $q\ge 4$, we have
\[
\frac d{dq}\Bigl[ q^{k-2}(q-1)+\frac 1{q^2}\Bigr]=(k-1)q^{k-2}-(k-2)q^{k-3}-2q^{-3}>0
\]
and
\[
\frac d{dq}\Bigl[ q^{k-2}(q-1)+\frac 2{q^2}\Bigr]=(k-1)q^{k-2}-(k-2)q^{k-3}-4q^{-3}\ge 0.
\]
Hence $f(q,k)$ is increasing with respect to $q$ for $q$ and $k$ in the above
range.
Thus, for $q\geq 4$ and $k\geq 1$, 
\[
f(q,k)\ge f(4,1)=\frac{819}{128}>6;
\]
for $q=3$ and $k\geq 2$,
\[
f(q,k)\ge f(3,2)=\frac{1520}{81}>6;
\]
for $q=2$ and $k\geq 4$,
\[
f(q,k)\ge f(2,4)=\frac{153}{8}>6.
\]
For $q=2$ and $k=3$, (\ref{3.4}) is verified directly:
\[
\frac{\genfrac{[}{]}{0pt}{0}{4} 3_2^2}{\dbinom{2^4}{2^3}} 
= \frac{5}{286}<\frac{1}{32}=\frac{1}{2^{2^3-3}}.
\]
\end{proof}

\begin{cor}\label{C3.3}
For $q\geq 4$, $n>k\geq 1$, or $q=3$, $n>k\geq 2$, or $q=2$, $n>k\geq 3$,
\begin{equation}\label{3.6}
\frac{\genfrac{[}{]}{0pt}{0} n k_q^2}{\dbinom{q^n}{q^k}} <
\frac{1}{q^{(n-k)(q^k-2k)+k}}.
\end{equation}
\end{cor}

\begin{proof}
Applying Lemma~\ref{L3.1} to the left hand side of (\ref{3.6}) $n-k-1$ times  and applying 
Lemma~\ref{L3.2} after that, we get 
\[
\begin{split}
\frac{\genfrac{[}{]}{0pt}{0} n k_q^2}{\dbinom{q^n}{q^k}} < \, &
\frac{1}{q^{(n-k-1)(q^k-2k)}}\cdot
\frac{\genfrac{[}{]}{0pt}{0} {k+1} k_q^2}{\dbinom{q^{k+1}}{q^k}}\cr
<\, &\frac{1}{q^{(n-k-1)(q^k-2k)}}
\cdot \frac{1}{q^{q^k-k}}\cr
=\, &\frac{1}{q^{(n-k)(q^k-2k)+k}}.
\end{split}
\]
\end{proof}
 
\begin{lem}\label{L3.4}
For $q\geq 4$, $k\geq 1$, or $q=3$, $k\geq 2$, or $q=2$, $k\geq 3$,
\begin{equation}\label{3.7} 
q^k-2k-2 \geq 0.
\end{equation}
\end{lem}

\begin{proof}
Let $g(q,k)$ denote the left hand side of (\ref{3.7}). 
For fixed $k\ge 1$, $g(q,k)$ is increasing with respect to $q$. Also note that 
\[
\frac{\partial g}{\partial k}=q^k\ln q-2>\frac{q^k-4}{2},
\] 
which is non-negative in the described range for $q$ and $k$. So, 
$g(q,k)$ is increasing with respect to $k$ in such range.
Thus, for $q\geq 4$ and $k\geq 1$
\[
g(q,k)\geq g(4,1)=0; 
\]
for $q=3$ and $k\geq 2$,
\[
g(q,k)\geq g(3,2)=3>0; 
\]
for $q=2$ and $k\geq 3$,
\[
g(q,k)\geq g(2,3)=0.  
\]
\end{proof}

\begin{proof}[Proof of Theorem~\ref{T2.3}] 
Let $q$ and $m$ satisfy one of the conditions (a) -- (c) in Theorem~\ref{T2.3} and let
$n>m$. By Corollary~\ref{C3.3} and Lemma~\ref{L3.4}, for $m\le k<n$, we have 
\[
q^{2(n-k)} \frac{\genfrac{[}{]}{0pt}{0} n k_q^2}{\dbinom{q^n}{q^k}}<
\frac{1}{q^{(n-k)(q^k-2k-2)+k}}\leq \frac{1}{q^k}.
\]
Therefore, 
\[
\sum_{k=m}^{n-1} q^{2(n-k)} \frac{\genfrac{[}{]}{0pt}{0} n k_q^2}{\dbinom{q^n}{q^k}} < 
\sum_{k=m}^{n-1}\frac{1}{q^k}<\sum_{k=1}^\infty\frac{1}{q^k}=\frac{1}{q-1}\leq 1,
\]
which proves the theorem.
\end{proof}

%%%%%%%%%%%%%%%%%%%% section 4 %%%%%%%%%%%%%%%%%%%%%%%%%%%%%%%%%%%%%%%%% 
\section{Examples of  $k$-Spectra for Small Values of $k,n,q$}

Recall that for given $k$, $n$, $q$,
\[
k\text{-spectrum}(n,q) = \{ k\text{-affinity}(f) :  f\in
\text{Per}({\mathbb F}_q^n) \}.
\]
Here we give  some examples. In all cases $a$ -- $b$ denotes all integers from
$a$ to $b$. Only a few of the spectra (as indicate below) are proved to be complete.
The other examples are spectra obtained by random searches. In such cases it is
possible but not certain  that some values are missing, hence the results are called
\emph{partial spectra}.
\vskip2mm

\noindent The full spectrum  for $k=1, n = 2, q = 3$:

\noindent$\{$0 -- 4, 6, 12$\}$

\medskip

\noindent The full spectrum  for $k=1, n = 2, q = 4$:

\noindent$\{$0 -- 6, 8, 12, 20$\}$

\medskip

\noindent A partial spectrum for $k=1, n = 2, q = 5$:

\medskip

\noindent $\{$0 -- 12, 14, 15, 20, 30$\}$

\medskip

\noindent A partial spectrum for $k=1, n = 2, q = 7$: 

\medskip

\noindent$\{$0 -- 2, 5 -- 7, 10 -- 12, 14 -- 18, 20 -- 30, 32, 35, 42, 
56$\}$

\medskip

\noindent A partial spectrum  for $k=1, n = 3, q = 3$:  
\medskip

\noindent$\{$0 -- 64, 66, 70, 71, 73, 75, 81, 93, 117$\}$

\medskip

\noindent The full spectrum for $k=2, n = 3, q = 2$:

\medskip

\noindent$\{$0, 2, 6, 14$\}$

\medskip

\noindent A partial spectrum  for $k=2, n = 3, q = 3$:
\medskip

\noindent $\{$0 -- 9, 11 -- 13, 15, 21, 39$\}$

\medskip

\noindent The full spectrum for $k=2, n = 4, q = 2$:
\medskip

\noindent
$\{$5 -- 20, 22, 24 -- 26, 28, 30, 32, 36, 38, 44, 48, 52, 56, 76, 84, 140$\}$

\medskip

\noindent A partial spectrum for $k=2, n = 5, q = 2$:

\medskip

\noindent $\{$0, 9 -- 416, 418, 420, 422, 424, 426 -- 428, 430 -- 432, 434, 
436 -- 440, 442, 444 -- 452, 454, 456 -- 462, 464, 466, 468 -- 472, 474, 476, 480, 
482, 484, 486, 488, 490, 492, 496, 500, 504, 506, 508, 512, 
514 -- 515, 517 -- 518, 
520, 526 -- 528, 530, 532, 536, 540, 548, 550, 552, 554, 556, 560, 564, 568, 
576, 600, 604 -- 605, 608, 618, 620, 640, 648, 664, 704, 706, 728, 732, 736,
792, 820, 960, 1240$\}$

\medskip

\noindent A partial spectrum for $k=2, n = 6, q = 2:$

\medskip

\noindent $\{$21, 28 -- 5132, 5134, 5136 -- 5140, 5142 -- 5148, 
5150 -- 5384,
5386, 5388, 5390, 5392, 5394,   5396 -- 5418, 5420, 
5422 -- 5446,  5448 -- 5464,  5466 -- 5468, 5470,  5472 -- 5480,
        5482,   5484 -- 5486, 5488, 5490, 5492, 5496, 5498, 5500, 5502,
        5504, 5506, 5508, 5510, 5512, 5514, 5516, 5520, 5522, 5524,
        5526, 5528, 5530, 5532, 5534, 5536, 5538, 5540, 5544, 5548,
          5550 -- 5552,   5554 -- 5556, 5558,   5560 -- 5562, 5564,
          5566 -- 5568,   5570 -- 5572, 5574, 5576, 5578, 5580 -- 5581,
        5584, 5586,   5588 -- 5606,   5608 -- 5610,   5612 -- 5634,
          5636 -- 5712,   5716 -- 5736, 5738, 5740 -- 5760, 5765 -- 5774,
          5776 -- 5784, 5786, 5788 -- 5790, 5792, 5794, 5796, 5820,
        5822, 5824, 5830, 5832, 5834, 5836, 5840,   5844, 5849,
          5855 -- 5857,   5860 -- 5861, 5868, 5874, 5876, 5878, 5880,
        5882, 5884, 5886, 5888, 5890, 5892, 5896, 5898, 5900, 5904,
        5908, 5912, 5936, 5940, 5944, 5948, 5952, 5960, 5974, 5976,
        5978, 5984, 5986, 5988, 5990, 5994, 6000, 6012,   6020, 6030,
        6032, 6034,   6045 -- 6058,   6070 -- 6084, 6086, 6088, 6090,
        6092, 6096,   6099 -- 6111,   6120 -- 6136, 6138, 6140, 6142,
        6144, 6146, 6148, 6152, 6160, 6178 -- 6180,  6182 -- 6184,
        6186,   6188 -- 6190, 6192, 6194, 6202, 6208, 6228, 6230, 6232,
        6234, 6236, 6238, 6240, 6242, 6244, 6248, 6256, 6260, 6264,
        6286, 6288, 6290, 6296, 6298, 6316, 6336, 6340, 6352, 6360,
        6384, 6444, 6448, 6452, 6460,  6475 -- 6480, 6482, 6484, 6488,
        6496,   6501 -- 6503,   6505 -- 6508,   6528 -- 6534, 6536, 6538,
        6540, 6544, 6557 -- 6560, 6586, 6588, 6590, 6592, 6594, 6596,
        6640, 6642, 6644, 6648, 6650, 6656, 6756, 6764, 6768, 6832,
        6955,   6957 -- 6958,  6982 -- 6984, 6986, 6988, 6992, 7036,
        7038, 7040, 7042, 7044, 7048, 7056, 7152, 7156, 7160, 7384,
        7461, 7490, 7492, 7552, 7994, 8052, 8056, 8176, 8556, 9176,
        10416$\}$

\medskip

\noindent {\bf Observations}:

\begin{enumerate}
\item 
From Theorem~\ref{T2.1}, $0\in k\text{-spectrum}(n,q)$ for $1\le k\le n-1$ unless 
$k=1$, $q=2$ or $k=2$, $q=2$ and $n$ is even.

\item 
Near the beginning in each example spectrum there is a long sequence of
consecutive values. For
$q=2$ and $k=2$, there seems to be a gap preceding the first non-zero affinity. For $q >2$ the
limited experimental data shows no such gaps. This suggests that for $q=2$, if
there is one
$2$-flat that is carried to a
$2$-flat then there are a certain number of other $2$-flats that must also be carried
to
$2$-flats.

\item 
Preceding the largest value in $k$-spectrum($n$,$q$) there appears to  a gap
of size
$2q^k \genfrac{[}{]}{0pt}{1}{n-1} k_q$. Note that this gap may be considered a threshold in the
sense that if $k\text{-affinity}(f) > q^{n-k}\genfrac{[}{]}{0pt}{1} n k_q - 2q^k \genfrac{[}{]}{0pt}{1}{n-1} k_q$, then $f$ takes all
$k$-flats to $k$-flats.
\end{enumerate}

%%%%%%%%%%%% section 4 1/2%%%%%%%%%%%%%%%%%%%%%%%%%%%%%%%%%%%%%
\section{$(n-1)\text{-spectrum}(n,2)$}

In this section, we will determine $(n-1)\text{-spectrum}(n,2)$, which is the set of all $(n-1)$-affinities
of permutations of ${\mathbb F}_2^n$. The standard dot product of $a,b\in{\mathbb F}_2^n$ is denoted by 
$\langle a,b\rangle$.
Every $(n-1)$-flat in ${\mathbb F}_2^n$ is uniquely of the form 
\[
H(a,\epsilon):=\{x\in {\mathbb F}_2^n: \langle a, x\rangle=\epsilon\}
\]
for some $a\in {\mathbb F}_2^n\setminus\{0\}$ and $\epsilon\in{\mathbb F}_2$. Let $f\in\text{Per}({\mathbb F}_2^n)$.
If for some $a\in {\mathbb F}_2^n\setminus\{0\}$ and some $\epsilon\in{\mathbb F}_2$, $f\bigl(H(a,\epsilon)\bigr)$ is an $(n-1)$-flat, say $f\bigl(H(a,\epsilon)\bigr)=H(b,\delta)$
for some $b\in{\mathbb F}_2^n\setminus\{0\}$ and $\delta\in{\mathbb F}_2$, we must have $f\bigl(H(a,1+\epsilon)\bigr)
=H(b,1+\delta)$. Therefore, for each such $a$ and $b$, there exists $\phi\in\text{Per}({\mathbb F}_2)$ such that
\begin{equation}\label{I5.1}
f\bigl(H(a,t)\bigr)=H\bigl(b,\phi(t)\bigr)\quad\text{for all}\ t\in {\mathbb F}_2.
\end{equation}

\begin{lem}\label{IL5.1}
Let $f\in\text{\rm Per}({\mathbb F}_2^n)$ and let 
\[
V_f=\{0\}\cup\bigl\{a\in {\mathbb F}_2^n\setminus\{0\}: f\bigl(H(a,0)\bigr)\ \text{is an $(n-1)$-flat}\bigr\}.
\]
Then $V_f$ is a subspace of ${\mathbb F}_2^n$
\end{lem}

\begin{proof}
For $a_1,a_2\in V_f$, we prove that $a_1+a_2\in V_f$. We may assume that $a_1\ne 0$, $a_2\ne 0$, and $a_1\ne a_2$.
By (\ref{I5.1}), there exist $b_i\in {\mathbb F}_2^n\setminus\{0\}$ and $\phi_i\in \text{Per}({\mathbb F}_2)$, 
$i=1,2$, such that
\[
f\bigl(H(a_i,t)\bigr)=H\bigl(b_i,\phi_i(t)\bigr)\quad\text{for all}\ t\in {\mathbb F}_2.
\]
Clearly, $b_1\ne b_2$. Since ${\mathbb F}_2$ has only two permutations, $t\mapsto t$ or $t\mapsto t+1$, we see that 
$\phi_1+\phi_2$ is a constant, say $\epsilon$. For any $x\in H(a_1+a_2,0)$, let $t=\langle a_1, x\rangle=
\langle a_2, x\rangle$. Then $x\in
H(a_i,t)$, hence $f(x)\in H\bigl(b_i,\phi_i(t)\bigr)$. It follows that
\[
\langle b_1+b_2,\, f(x)\rangle=\phi_1(t)+\phi_2(t)=\epsilon,
\]
i.e., $f(x)\in H(b_1+b_2,\epsilon)$. Thus we have proved that $f\bigl(H(a_1+a_2,0)\bigr)=H(b_1+b_2,\epsilon)$,
which implies that $a_1+a_2\in V_f$.
\end{proof}

\begin{thm}\label{IT5.2}
Let $n>2$. Then
\[
(n-1)\text{\rm -spectrum}(n,2)=\{2^i-2: 1\le i\le n+1\}.
\]
\end{thm}

\begin{proof}
For each $f\in \text{Per}({\mathbb F}_2^n)$, by Lemma~\ref{IL5.1}, we have 
\[
(n-1)\text{-affinity}(f)=2|V_f\setminus\{0\}|=2^{\dim V_f+1}-2\in \{2^i-2: 1\le i\le n+1\}.
\]
It remains to show that for each $1\le i\le n+1$, there exists $f\in \text{Per}({\mathbb F}_2^n)$ with
\[
(n-1)\text{-affinity}(f)=2^i-2.
\]
We prove this claim by induction on $n$. For $n=3$, the claim was established by computer as mentioned in 
Section 4. Assume $n>3$. If $i=1$, the claim follows from Theorem~\ref{T2.1}. Thus we will assume $2\le i\le n+1$. By the induction hypothesis, there
exists $g\in \text{Per}({\mathbb F}_2^{n-1})$ such that $(n-2)\text{-affinity}(g)=2^{i-1}-2$. Define $f\in
\text{Per}({\mathbb F}_2^n)$ by $f(c,x)=(c,g(x))$, $c\in{\mathbb F}_2$, $x\in{\mathbb F}_2^{n-1}$. Clearly,
$\{i\}\times{\mathbb F}_2^{n-1}$, $i=0,1$, are mapped into flats by $f$. Let $X\subset{\mathbb F}_2^n$ be any
$(n-1)$-flat other than $\{i\}\times{\mathbb F}_2^{n-1}$, $i=0,1$, such that $f(X)$ is a flat. Write
\[
X\cap\bigl(  \{i\}\times{\mathbb F}_2^{n-1} \bigr)=\{i\}\times U_i,\quad i=0,1.,
\]
where $U_i\subset{\mathbb F}_2^{n-1}$ is an $(n-2)$-flat and $U_0=U_1$ or ${\mathbb F}_2^{n-1}\setminus U_1$.
Then
\begin{equation}\label{I5.2}
X=\bigl(\{0\}\times U_0\bigr)\cup\bigl(\{1\}\times U_1\bigr).
\end{equation}
Since 
\[
\{i\}\times g((U_i)=f\Bigl(X\cap\bigl(  \{i\}\times{\mathbb F}_2^{n-1} \bigr)\Bigr)=f(X)\cap\bigl(  \{i\}\times{\mathbb F}_2^{n-1} \bigr)
\]
is an $(n-2)$-flat, $g(U_i)$ is an $(n-2)$-flat. On the other hand, given any $(n-2)$-flats $U_0,U_1\subset
{\mathbb F}_2^{n-1}$ such that $U_0= U_1$ or ${\mathbb F}_2^{n-1}\setminus U_1$ and $g(U_i)$, $i=0,1$, are
$(n-2)$-flats, both $X$ (in (\ref{I5.2})) and $f(X)$ are $(n-1)$-flats in ${\mathbb F}_2^n$. Therefore,
\[
\begin{split}
(n-1)\text{-affinity}(f)\,&=2+2\cdot\bigl((n-2)\text{-affinity}(g)\bigr)\cr
&=2+2(2^{i-1}-2)\cr
&=2^i-2.\cr
\end{split}
\]
The proof is now complete.
\end{proof}

%%%%%%%%%%%%% section 5 %%%%%%%%%%%%%%%%%%%%%%%%%%%%%%%%%%%%%%%%
\section{Threshold Conjecture for $k$-Affinity}

In many cases it appears that the next to largest $k$-affinity is the $k$-affinity of
a transposition. We calculate this value in the following lemma. In this case it is
more convenient to compute the $k$-coaffinity.

\begin{lem}\label{L5.1}
Let $f\in\text{\rm Per}({\mathbb F}_q^n)$ be any transposition. 
Then 
\[
k\text{\rm-coaffinity}(f) = 
 2q^k\genfrac{[}{]}{0pt}{0}{n-1} k_q
\]
and
\[ 
k\text{\rm-affinity}(f) = 
 \left( \frac{(q^{n-k}-2)(q^n-1)}{q^k-1} +2 \right)\genfrac{[}{]}{0pt}{0}{n-1}{k-1}_q.
\]
\end{lem}

\begin{proof}
Assume that $f$ 
interchanges $x$ and $y$, where $x,y\in{\mathbb F}_q^n$ are distinct.
If $U$ is a $k$-flat, then $f(U)$ is not a $k$-flat if and only if $U$ contains
exactly one of $x$ and $y$. 

The number of $k$-flats in ${\mathbb F}_q^n$ containing a fixed point is
$\genfrac{[}{]}{0pt}{1}{n} k_q$; the number of $k$-flats in ${\mathbb F}_q^n$ containing two fixed points is
$\genfrac{[}{]}{0pt}{1}{n-1} {k-1}_q$.
Hence
\[
k\text{-coaffinity}(f)=
 2\left( \genfrac{[}{]}{0pt}{0}{n} k_q - \genfrac{[}{]}{0pt}{0}{n-1}{k-1}_q\right) =
2q^k \genfrac{[}{]}{0pt}{0}{n-1} k_q.
\]
It follows that
\[
\begin{split}
k\text{-affinity}(f)\,& = 
q^{n-k}\genfrac{[}{]}{0pt}{0}{n} k_q - 2q^k\genfrac{[}{]}{0pt}{0}{n-1} k_q\cr 
&= \left( \frac{(q^{n-k}-2)(q^n-1)}{q^k-1} + 2 \right)\genfrac{[}{]}{0pt}{0}{n-1}{k-1}_q.
\end{split}
\]
\end{proof}

It is natural to call the next largest $k$-affinity a threshold for affinity and we
make the following conjecture:

\begin{conj}[Threshold Conjecture]\label{CJ4.2}
Let $1\le k\le n-1$ and $f \in \text{\rm Per}({\mathbb F}_q^n)$.
If
\[
k\text{\rm -coaffinity}(f)< 2q^k\genfrac{[}{]}{0pt}{0}{n-1} k_q,
\]
then $k\text{\rm -coaffinity}(f)=0$, i.e., $f\in{\rm AGL}(n,{\mathbb F}_q)$. Equivalently, if
\[
k\text{\rm -affinity}(f)> \left( \frac{(q^{n-k}-2)(q^n-1)}{q^k-1} + 2 \right)\genfrac{[}{]}{0pt}{0}{n-1}{k-1}_q,
\]
then $k\text{\rm -affinity}(f)=q^{n-k}\genfrac{[}{]}{0pt}{1}{n} k_q$. 
That is, the next to largest $k$-affinity is that of a transposition. 
\end{conj}

This conjecture is supported by the examples in Section 4 and the result in Section 5.
More importantly it is supported by the proof for 
 $q=2$, $k = 2$, $n > 2$  in Section 7,  and the proof for 
$q>2$, $k=1$, $n>1$  in Section 8.

%%%%%%%%%%%%%%%%%%%%%%%% section 6 %%%%%%%%%%%%%%%%%%%%%%%%%%%%%%%%%%%%%
\section{Proof of the Threshold Conjecture for $k=2$, $q=2$}

Recall that a 2-flat in ${\mathbb F}_2^n$ is simply a 4-element subset $\{x_1,x_2,x_3,x_4\}$
such that $x_1+x_2+x_3+x_4=0$. For $f\in\text{Per}({\mathbb F}_2^n)$ and a 2-flat $X\subset
{\mathbb F}_2^n$, $f(X)$ is a 2-flat if and only if $f$ is affine on $X$.

For the proof in this section, the reader's familiarity with the Fourier transformation of boolean functions will be helpful. We first introduce the necessary notation.
The set of all functions from ${\mathbb F}_2^n$ to
${\mathbb F}_2$ is denoted by ${\mathcal P}_n$. 
Every function in ${\mathcal P}_n$ is uniquely represented by a polynomial in
${\mathbb F}_2[X_1,\dots, X_n]$ whose degree in each $X_i$ is at most 1. Namely,
\[
{\mathcal P}_n={\mathbb F}_2[X_1,\dots,X_n]/\langle X_1^2-X_1,\dots, X_n^2-X_n
\rangle.
\]
For each $g\in{\mathcal P}_n$, put $|g|=|g^{-1}(1)|$.
The Fourier transform of $g\in{\mathcal P}_n$ is the function $\hat g:{\mathbb F}_2^n\rightarrow{\mathbb C}$ defined
by
\[
\hat g(a)=\sum_{x\in{\mathbb F}_2^n}(-1)^{g(x)+\langle a, x\rangle},\quad a\in{\mathbb F}_2^n,
\]
where $\langle a, x\rangle$ is the standard dot product in ${\mathbb F}_2^n$. Clearly,
\begin{equation}\label{E6.0}
\hat g(a)=2^n-2|g+\langle a,\, \cdot\rangle|.
\end{equation}
Note that for $n\ge 2$, $|g+\langle a,\, \cdot\rangle|\equiv |g|\pmod 2$, hence
\[
\hat g(a)\equiv 2|g|\pmod 4.
\]
It is well known (also straightforward to prove) that
\begin{equation}\label{E6.1}
\sum_{a\in{\mathbb F}_2^n}\bigl(\hat g(a)\bigr)^2=2^{2n}
\end{equation}
and 
\begin{equation}\label{E6.2}
\sum_{a\in{\mathbb F}_2^n}\bigl(\hat g(a)\bigr)^4=2^n\sum_{a\in{\mathbb F}_2^n}\Bigl[
\sum_{x\in{\mathbb F}_2^n} (-1)^{g(x+a)+g(x)}\Bigr]^2.
\end{equation}
(Equation (\ref{E6.1}) is the Parseval identity; equation (\ref{E6.2}) is a relation between the Fourier transform
and the convolution of the function. Cf. \cite{Lan91}.) If $A\le \bigl(\hat g(a)\bigr)^2\le B$ for all 
$a\in{\mathbb F}_2^n$, from (\ref{E6.1}), we have
\[
\begin{split}
2^n\Bigl(\frac{B-A}2\Bigr)^2\,&\ge \sum_{a\in{\mathbb F}_2^n}\Bigl[\bigl(\hat g(a)\bigr)^2-\frac{A+B}2\Bigr]^2\cr
&=\sum_{a\in{\mathbb F}_2^n}\bigl(\hat g(a)\bigr)^4-(A+B)\sum_{a\in{\mathbb F}_2^n}\bigl(\hat g(a)\bigr)^2+2^n
\Bigl(\frac{A+B}2\Bigr)^2\cr
&=\sum_{a\in{\mathbb F}_2^n}\bigl(\hat g(a)\bigr)^4-2^{2n}(A+B)+2^n\Bigl(\frac{A+B}2\Bigr)^2.
\end{split}
\]
Thus 
\[
\sum_{a\in{\mathbb F}_2^n}\bigl(\hat g(a)\bigr)^4\le 2^{2n}(A+B)-2^nAB.
\]
Combining the above with (\ref{E6.2}), we have
\begin{equation}\label{E6.4}
\sum_{a\in{\mathbb F}_2^n}\Bigl[
\sum_{x\in{\mathbb F}_2^n} (-1)^{g(x+a)+g(x)}\Bigr]^2\le 2^n(A+B)-AB.
\end{equation}
The equality in (\ref{E6.4}) holds if and only if $\bigl(\hat g(a)\bigr)^2=A$ or $B$ for all $a\in{\mathbb F}_2^n$.

\begin{lem}\label{L6.1}
Let $g\in{\mathcal P}_n$ with $\deg g\ge 2$. Then
\begin{equation}\label{6.1}
\sum_{a\in{\mathbb F}_2^n}\Bigl[
\sum_{x\in{\mathbb F}_2^n}(-1)^{g(x+a)+g(x)}\Bigr]^2\le
2^{2n}+(2^n-1)(2^n-4)^2.
\end{equation}
The equality holds if and only if $|g+h|=1$ 
for some $h\in{\mathcal P}_n$ with $\deg h\le 1$.
\end{lem}

\begin{proof}
Since $\deg g\ge 2$, $\hat g(a)\ne \pm 2^n$ for all $a\in{\mathbb F}_2^n$.

{\bf Case 1}. $|g|$ is even. By (\ref{E6.0}), $0\le |\hat g(a)|\le 2^n-4$ for  all $a\in{\mathbb F}_2^n$.
Hence
\[
0\le \big(\hat g(a)\bigr) ^2\le (2^n-4)^2\quad\text{for all}\ a\in{\mathbb F}_2^n.
\]
By (\ref{E6.4}),
\[
\sum_{a\in{\mathbb F}_2^n}\Bigl[
\sum_{x\in{\mathbb F}_2^n}(-1)^{g(x+a)+g(x)}\Bigr]^2\le
2^n(2^n-4)^2<2^{2n}+(2^n-1)(2^n-4)^2.
\]

{\bf Case 2}. $|g|$ is odd, By (\ref{E6.0}), $2\le|\hat g(a)|\le 2^n-2$ for  all $a\in{\mathbb F}_2^n$.  
Hence
\[
2^2\le \big(\hat g(a)\bigr) ^2\le (2^n-2)^2\quad\text{for all}\ a\in{\mathbb F}_2^n.
\]
By (\ref{E6.4}),
\[
\begin{split}
\sum_{a\in{\mathbb F}_2^n}\Bigl[
\sum_{x\in{\mathbb F}_2^n}(-1)^{g(x+a)+g(x)}\Bigr]^2\,&\le 2^n\bigl[2^2+(2^n-2)^2\bigr]-2^2(2^n-2)^2\cr
&=2^{2n}+(2^n-1)(2^n-4)^2.
\end{split}
\]
The equality holds if and only if $\bigl(\hat g(a)\bigr)^2=2^2$ or $(2^n-2)^2$ for all $a\in{\mathbb F}_2^n$.

First assume that the equality in (\ref{6.1}) holds.
Since $\sum_{a\in{\mathbb F}_2^n}\bigl(\hat g(a)\bigr)^2=2^{2n}$, for at least one $a\in{\mathbb F}_2^n$,
$\bigl(\hat g(a)\bigr)^2=(2^n-2)^2$.  By (\ref{E6.0}),
\[
\bigl(2^n-2|g+\langle a,\,\cdot\rangle|\bigr)^2=(2^n-2)^2.
\]
Thus $|g+\langle a,\,\cdot\rangle|=1$ or $2^n-1$. Let
\[
h=
\begin{cases}
\langle a,\,\cdot\rangle&\text{if}\ |g+\langle a,\,\cdot\rangle|=1,\cr
\langle a,\,\cdot\rangle+1&\text{if}\ |g+\langle a,\,\cdot\rangle|=2^n-1.\cr
\end{cases}
\]
Then $|g+h|=1$, as claimed.

Now assume that $|g+h|=1$ for some $h\in{\mathcal P}_n$ with $\deg h\le 1$. Then for every $a\in{\mathbb F}_2^n$,
\[
|g+\langle a,\,\cdot\rangle|=2^{n-1}\pm 1,\ \text{or}\ 1, \text{or}\ 2^n-1.
\]
It follows from (\ref{E6.0}) that $\hat g(a)=\pm 2$ or $\pm(2^n-2)$. Hence $\bigl(\hat g(a)\bigr)^2=2^2$ or $(2^n-2)^2$ for all $a\in{\mathbb F}_2^n$. Therefore the equality holds in (\ref{6.1}).
\end{proof}

\begin{thm}\label{T6.2}
Let $n\ge 3$ and $f\in \text{\rm Per}({\mathbb F}_2^n)\setminus{\rm AGL}(n,{\mathbb F}_2)$. Then
\[
2\text{\rm -coaffinity}(f)\ge \frac 83 (2^{n-1}-1)(2^{n-2}-1).
\]
The equality holds if and only if $f\in{\rm AGL}(n,{\mathbb F}_2)\circ \tau\circ
{\rm AGL}(n,{\mathbb F}_2)$ where $\tau\in {\rm Per}({\mathbb F}_2^n)$ is any transposition.
\end{thm}

\begin{cor}\label{C6.3}
The Threshold Conjecture (Conjecture~\ref{CJ4.2}) holds for $q=2$, $k=2$, $n>2$.
\end{cor}

\begin{proof}[Proof of Theorem~\ref{T6.2}]
{\bf Case 1}. $f(x+d)+f(x)=$ constant for some $d\in{\mathbb F}_2^n\setminus\{0\}$.
Without loss of generality, assume $d=(1,0,\dots,0)$. Write $f=(f_1,\dots,f_n)$ where
$f_i:{\mathbb F}_2^n\rightarrow {\mathbb F}_2$. Then
\[
f_i=c_iX_1+g_i(X_2,\dots,X_n)
\]
for all $1\le i\le n$,
where $c_i\in{\mathbb F}_2$ and $g_i\in{\mathcal P}_{n-1}$. Note that $(c_1,\dots,c_n)\ne 0$ since
otherwise, $f$ is independent of $X_1$ and cannot be a permutation of ${\mathbb F}_2^n$. By
composing a suitable element of ${\rm GL}(n,{\mathbb F}_2)$ to the left of $f$, we may assume 
$(c_1,\dots,c_n)=(1,0,\dots,0)$. Thus
\[
f=\bigl(X_1+g_1(X_2,\dots,X_n),\ g_2(X_2,\dots,X_n),\ \dots,\ g_n(X_2,\dots,X_n)\bigr).
\]
It follows that $g=(g_2,\dots,g_n)$ is a permutation of ${\mathbb F}_2^{n-1}$.

{\bf Case 1.1}. $g\notin {\rm AGL}(n-1,{\mathbb F}_2)$. Using induction, we may assume
\[
2\text{-coaffinity}(g)\ge \frac 83 (2^{n-2}-1)(2^{n-3}-1).
\]
Note that if $g$ is not affine on a 2-flat $\{y_1,y_2,y_3,y_4\}$ in ${\mathbb F}_2^{n-1}$,
$f$ is not affine on the 2-flat $\{(x_i,y_i):i\le i\le 4\}$ where $x_1+\cdots+x_4=0$. Hence
\[
2\text{-coaffinity}(f)\ge 8\cdot\bigl(2\text{-coaffinity}(g)\bigr)>\frac 83 (2^{n-1}-1)(2^{n-2}-1).
\]

{\bf Case 1.2}. $g\in {\rm AGL}(n-1,{\mathbb F}_2)$. Then $\deg g_1\ge 2$. We may assume $g=
{\rm id}$. In this case, the 2-flats on which $f$ is not affine are precisely
\[
\{(x_i,y_i):1\le i\le 4\}
\]
where $\{y_1,\dots, y_4\}$ is a 2-flat in ${\mathbb F}_2^{n-1}$ such that $\sum_{i=1}^4 g_1(y_i)=1$
and $x_1,\dots,x_4\in{\mathbb F}_2$ with $\sum_{i=1}^4x_i=0$. Define
\[
\begin{array}{crcl}
G:&({\mathbb F}_2^{n-1})^3&\longrightarrow &{\mathbb F}_2\cr
& (y,a,b)& \longmapsto & g_1(y)+g_1(y+a)+g_1(y+b)+g_1(y+a+b)\cr
\end{array}
\]
Then 
\[
2\text{-coaffinity}(f)=\frac 8{4!} |G|=\frac 13 |G|.
\]
We have
\[
\begin{split}
|G|\,&=\frac 12\Bigl[ 2^{3(n-1)}-\sum_{y,a,b\in{\mathbb F}_2^{n-1}}(-1)^{G(y,a,b)}\Bigr]\cr
&=\frac 12\Bigl[ 2^{3(n-1)}-\sum_{a\in{\mathbb F}_2^{n-1}}
\sum_{y,b\in{\mathbb F}_2^{n-1}}(-1)^{g_1(y)+g_1(y+a)+g_1(y+b)+g_1(y+a+b)}\Bigr]\cr
&=\frac 12\Bigl[ 2^{3(n-1)}-\sum_{a\in{\mathbb F}_2^{n-1}}\Bigl(
\sum_{y\in{\mathbb F}_2^{n-1}}(-1)^{g_1(y)+g_1(y+a)}\Bigr)^2\Bigr]\cr
& \ge \frac 12\Bigl[ 2^{3(n-1)}-\bigl[2^{2(n-1)}+(2^{n-1}-1)(2^{n-1}-4)^2\bigr]\Bigr]\quad
\text{(by Lemma~\ref{L6.1})}\cr
&=2^3(2^{n-1}-1)(2^{n-2}-1).\cr
\end{split}
\]
Therefore,
\[
2\text{-coaffinity}(f)=\frac13 |G|\ge \frac 83(2^{n-1}-1)(2^{n-2}-1).
\]

If the equality holds in the above, then the equality in (\ref{6.1}) holds with $g_1$ in place of $g$. By Lemma~\ref{L6.1}, there exists $h\in{\mathcal P}_{n-1}$ such that 
$|g_1+h|=1$. Using a linear transformation, we may replace $g_1$ with $g_1+h$. Thus
we may assume $|g_1|=1$. Then clearly,
\[
f=\bigl(X_1+g_1(X_2,\dots,X_n),\ X_2,\ \dots,\ X_n\bigr)
\]
is a transposition.

{\bf Case 2}. $f(x+d)+f(x)\ne$ constant for all $d\in{\mathbb F}_2^n\setminus\{0\}$. For each 
$d\in{\mathbb F}_2^n\setminus\{0\}$, let
\[
\Delta(d)=\bigl\{\{x,x+d\}:x\in{\mathbb F}_2^n\bigr\}\subset\binom{{\mathbb F}_2^n}2,
\]
where $\binom{{\mathbb F}_2^n}2$ denotes the set of all 2-element subsets of ${\mathbb F}_2^n$. 
Denote $\{f(X):X\in\Delta(d)\}$ by $f(\Delta(d))$ (an abuse of notation for the convenience).
By the assumption, $f\bigl( \Delta(d) \bigr)\not\subset \Delta(c)$ for every 
$c\in{\mathbb F}_2^n\setminus\{0\}$. Since the subsets in $\Delta(d)$ and $\Delta(c)$ form partitions of ${\mathbb F}_2^n$, we have 
\begin{equation}\label{6.10}
\bigl| f\bigl( \Delta(d)\bigr) \cap \Delta(c)
\bigr|\le 2^{n-1}-2.
\end{equation}

We claim that we can partition ${\mathbb F}_2^n\setminus\{0\}$ into $A$ and $B$ such that
\[
\Bigl|f\bigl( \Delta(d) \bigr)\cap\Bigl(\bigcup_{a\in A}\Delta(a)\Bigr)\Bigr|\ge 2,
\]
\[
\Bigl|f\bigl( \Delta(d) \bigr)\cap\Bigl(\bigcup_{b\in B}\Delta(b)\Bigr)\Bigr|\ge 2.
\]
Note that $\Delta(a)$, $a\in{\mathbb F}_2^n\setminus\{0\}$ form a partition of $\binom
{{\mathbb F}_2^n} 2$. If $\bigl|f\bigl( \Delta(d) \bigr)\cap\Delta(a)\bigr|\le 1$
for all $a\in{\mathbb F}_2^n\setminus\{0\}$, choose $a_1,a_2\in{\mathbb F}_2^n\setminus\{0\}$
distinct such that $\bigl|f\bigl( \Delta(d) \bigr)\cap\Delta(a_i)\bigr|= 1$, $i=1,2$.
Then $A=\{a_1,a_2\}$, $B={\mathbb F}_2^n\setminus\{0, a_1, a_2\}$ have the desired
property. If $\bigl|f\bigl( \Delta(d) \bigr)\cap\Delta(a)\bigr|\ge 2$ for some $a\in
{\mathbb F}_2^n\setminus\{0\}$, let $A=\{a\}$ and $B={\mathbb F}_2^n\setminus\{0, a\}$.
By (\ref{6.10}), we have
\[
\Bigl|f\bigl( \Delta(d) \bigr)\cap\Bigl(\bigcup_{b\in B}\Delta(b)\Bigr)\Bigr|=
2^{n-1}-\bigl|f\bigl( \Delta(d) \bigr)\cap\Delta(a)\bigr|\ge 2.
\]
Hence $A$ and $B$ also have the desired property.

Therefore, among the 2-flats which are a union of two elements in $\Delta(d)$, there are at least
\[
2\cdot(2^{n-1}-2)
\]
on which $f$ is not affine. Since this statement is true for all $d\in{\mathbb F}_2^n\setminus\{0\}$, it follows that 
\[
2\text{-coaffinity}(f)\ge\frac{2\cdot (2^{n-1}-2)\cdot (2^n-1)}3>\frac 83(2^{n-1}-1)(2^{n-2}-1).
\]
\end{proof}

For the remainder of this section, we compute the number of permutations $f\in\text{Per}
({\mathbb F}_2^n)$ with $2\text{-coaffinity}(f)=2^2
\genfrac{[}{]}{0pt}{1}{n-1} 2_2=\frac 83(2^{n-1}-1)(2^{n-2}-1)$.

\begin{lem}\label{L6.4}
Let $n\ge 3$ and choose $a\in{\mathbb F}_2^n\setminus\{0\}$. Let $\tau\in {\rm Per}({\mathbb F}_2^n)$ be the transposition 
which permutes $0$ and $a$ and let $f\in{\rm AGL}(n,{\mathbb F}_2)$.
\begin{enumerate}
\item[(i)] If $n\ge 4$, then $\tau\circ f\circ\tau\in{\rm AGL}(n,{\mathbb F}_2)$ if and only if
$f(\{0,a\})=\{0,a\}$.

\item[(ii)] If $n=3$, then $\tau\circ f\circ\tau\in{\rm AGL}(n,{\mathbb F}_2)$ if and only if
$f(a)+f(0)=a$.
\end{enumerate}
\end{lem}

\begin{proof}
(i) ($\Leftarrow$) In fact, $\tau\circ f\circ\tau=f$.

($\Rightarrow$) Assume to the contrary that $f(\{0,a\})\ne\{0,a\}$. Without loss of generality, assume $f^{-1}(0)
\notin\{0,a\}$. Since $n\ge 4$, there is a 2-flat $A$ in ${\mathbb F}_2^n$ which contains $f^{-1}(0)$ but does
contain $0, a$ and $f^{-1}(a)$. Write $A=\{f^{-1}(0), b_2,b_3,b_4\}$ where $b_i\notin\{0,a\}$, $f(b_i)\notin
\{0,a\}$ and $f(b_2)+f(b_3)+f(b_4)=0$. We then have
\[
(\tau\circ f\circ\tau)(A)=\tau\bigl(f(A)\bigr)=\tau\bigl(\{0,f(b_2),f(b_3),f(b_4)\}\bigr)=\{a,f(b_2),f(b_3),f(b_4)\},
\]
which is not a 2-flat. This is a contradiction.

(ii) ($\Leftarrow$) Without loss of generality, assume $a=(1,0,0)$. Then
\[
\tau(x_1,x_2,x_3)=(x_1+g(x_2,x_3),\, x_2,\, x_3)
\]
where
\begin{equation}\label{6.11}
g(x_2,x_3)=
\begin{cases}
1&\text{if}\ (x_2,x_3)=0,\cr
0&\text{if}\ (x_2,x_3)\ne 0.\cr
\end{cases}
\end{equation}
Since $f(a)+f(0)=a$, we have
\[
f(x)=xA+b,
\]
where $b\in{\mathbb F}_2^3$, $A\in{\rm GL}(3,{\mathbb F}_2)$ and $aA=a$, i.e.,
\[
A=\left[
\begin{matrix}
1&0\cr
c&B\cr
\end{matrix}\right],\quad B\in{\rm GL}(2,{\mathbb F}_2),\ c\in{\mathbb F}_2^2.
\]
Let $b=(b',b'')$ where $b'\in{\mathbb F}_2$ and $b''\in{\mathbb F}_2^2$. Then
\[
\begin{split}
&(\tau\circ f\circ\tau)(x_1,x_2,x_3)\cr
=\,&(\tau\circ f)(x_1+g(x_2,x_2),\, x_2,\, x_3)\cr
=\,&\tau\Bigl((x_1+g(x_2,x_2),\, x_2,\, x_3)A+b\Bigr)\cr
=\,&\tau\Bigl(x_1+g(x_2,x_3)+(x_2,x_3)c+b',\ (x_2,x_3)B+b''\Bigr)\cr
=\,&\Bigl(x_1+g(x_2,x_3)+(x_2,x_3)c+b'+g\bigl((x_2,x_3)B+b''\bigr),\ (x_2,x_3)B+b''\Bigr).
\end{split}
\]
By (\ref{6.11}), $g\bigl((x_2,x_3)B+b''\bigr)=g\bigl((x_2,x_3)+b''B^{-1}\bigr)$. Thus
\[
g(x_2,x_3)+g\bigl((x_2,x_3)B+b''\bigr)=g(x_2,x_3)+g\bigl((x_2,x_3)+b''B^{-1}\bigr)
\]
has degree $\le 1$ since $\deg g\le 2$. Therefore $\tau\circ f\circ\tau\in{\rm AGL}(3,{\mathbb F}_2)$.

($\Rightarrow$) Assume to the contrary that $f(a)+f(0)\ne a$. 
Then $f\bigl(\{0,a\}\bigr)\ne\{0,a\}$. Without loss of generality, we may assume $f^{-1}(0)\notin\{0,a\}$.
By the proof of ($\Rightarrow$) of (i), it suffices
to show that there is a 2-flat $A$ in ${\mathbb F}_2^3$ which contains $f^{-1}(0)$ but not $0, a$ and $f^{-1}(a)$. 
By the assumption, the set $\{0,a,f^{-1}(0),f^{-1}(a)\}$ is not 2-flat of ${\mathbb F}_2^3$, hence it either has only
three distinct elements or is a frame of ${\mathbb F}_2^3$. 
(A set of $k+1$ elements in an affine space is called frame if their affine span is a $k$-flat.) 
To see the existence of a desirable 2-flat $A$, first note
that $\bigl\{(0,0,0),(1,0,0),(0,1,0),(0,0,1)\bigr\}$ is a frame of ${\mathbb F}_2^3$ and $\bigl\{(x_1,x_2,x_3):
x_1+x_2+x_3=0\bigr\}$ is a 2-flat which contains exactly one of the elements in the frame. Using a suitable 
affine transformation, we see that for any frame $\{\epsilon_0,\dots,\epsilon_3\}$ of ${\mathbb F}_2^3$, there is a 2-flat 
$A$ such that $A\cap\{\epsilon_0,\dots,\epsilon_3\}=\{\epsilon_0\}.$
\end{proof}

\begin{cor}\label{C6.5}
In the notation of Lemma~\ref{L6.4}, we have
\begin{equation}\label{6.12}
\bigl|\bigl(\tau\circ{\rm AGL}(n,{\mathbb F}_2)\circ\tau\bigr)\cap{\rm AGL}(n,{\mathbb F}_2)\bigr|
=
\begin{cases}
2^n|{\rm GL}(n-1,{\mathbb F}_2)|&\text{if}\ n\ge 4,\cr
2^5|{\rm GL}(2,{\mathbb F}_2)|&\text{if}\ n=3.\cr
\end{cases}
\end{equation}
\end{cor}

\begin{proof}
Let $f\in{\rm AGL}(n,{\mathbb F}_2)$ be given by
\[
f(x)=xA+b,
\]
where $A\in{\rm GL}(n,{\mathbb F}_2)$ and $b\in{\mathbb F}_2^n$. By Lemma~\ref{L6.4}, when $n\ge 4$, $f\in 
\tau\circ{\rm AGL}(n,{\mathbb F}_2)\circ\tau$ if and only if $aA=a$ and $b=0$ or $a$;  when $n=3$, $f\in 
\tau\circ{\rm AGL}(n,{\mathbb F}_2)\circ\tau$ if and only if $aA=a$. Equation (\ref{6.12}) follows immediately.
\end{proof}

\begin{cor}\label{C6.6}
Assume $n\ge 3$. The number of permutations $f\in{\rm Per}({\mathbb F}_2^n)$ 
with $2\text{\rm -coaffinity}(f)=2^2
\genfrac{[}{]}{0pt}{1}{n-1} 2_2$ is given by 
\[
\begin{cases}
2^{\frac 12(n^2+3n-2)}(2^n-1)^2\prod_{j=1}^{n-1}(2^j-1)&\text{if}\ n\ge 4,\cr
2^6\cdot 3\cdot 7^2&\text{if}\ n=3.\cr
\end{cases}
\]
\end{cor}

\begin{proof}
Let $\tau\in {\rm Per}({\mathbb F}_2^n)$ be any transposition. By Theorem~\ref{T6.2}, 
the number of $f\in{\rm Per}({\mathbb F}_2^n)$ with $2\text{-coaffinity}(f)=2^2
\genfrac{[}{]}{0pt}{1}{n-1} 2_2$ is
\[
|{\rm AGL}(n,{\mathbb F}_2)\circ\tau\circ {\rm AGL}(n,{\mathbb F}_2)|
=\frac{|{\rm AGL}(n,{\mathbb F}_2)|^2}{\bigl|\bigl(\tau\circ {\rm AGL}(n,{\mathbb F}_2)\circ\tau\bigr)
\cap{\rm AGL}(n,{\mathbb F}_2)\bigr|}.
\]
The result follows immediately from the above corollary.
\end{proof}

%%%%%%%%%%%%%%%%%% section 7 %%%%%%%%%%%%%%%%%%%%%%%%%%%%%%%%%%%%%%%%%%%
  
\section{Proof of the Threshold Conjecture for $k=1$, $q > 2$}

The Threshold Conjecture for $k=1$, $q>2$ is proved in two steps: first $n=2$ then $n\ge 3$.
For $n=2$, the cases $q=3$ and $q\ge 4$ require different treatments. The group of all invertible semi-affine transformations of ${\mathbb F}_q^n$ is denoted by ${\rm A\Gamma L}
(n,{\mathbb F}_q)$. For any two distinct points $x,y\in{\mathbb F}_q^n$, $\overline{xy}$
denotes the line through $x$ and $y$ in ${\mathbb F}_q^n$.

\begin{thm}\label{T7.1}
Let $q\ge 4$ and $f\in {\rm Per}({\mathbb F}_q^2)\setminus{\rm A\Gamma L}(2,{\mathbb F}_q)$.
Then
\[
1\text{\rm -coaffinity}(f)\ge 2q\,\genfrac{[}{]}{0pt}{0}1 1_q=2q.
\]
\end{thm}

\begin{proof}
Among the $q+1$ parallel classes of lines in ${\mathbb F}_q^2$, we first assume that at most one parallel class has the property that all lines
in the class are mapped to lines by $f$. In each of the remaining $q$
parallel classes, there are at least 2 lines which are not mapped to lines by $f$. (Since the lines in a parallel class form a partition of
${\mathbb F}_q^2$, it cannot be the case that exactly one line in a parallel
class is not mapped to a line.) Therefore $1\text{-coaffinity}(f)\ge 2q$.

Now assume that there are two parallel classes of lines in ${\mathbb F}_q^2$
such that all lines in the two parallel classes are mapped to lines by $f$.
By composing suitable linear transformations to both sides of $f$, we may assume that $f$ maps all horizontal lines to horizontal lines and all
vertical lines to vertical lines. (A horizontal line in 
${\mathbb F}_q^2$ is a line with direction vector $(1,0)$; a vertical line in 
${\mathbb F}_q^2$ is a line with direction vector $(0,1)$.)

Assume that for every $z\in{\mathbb F}_q^2$, there is a line through $z$ which is not mapped to a line by $f$. Then there are at least two lines through $z$ which are not mapped to lines.
Since each line contains $q$ points, we have
\[
1\text{-coaffinity}(f)\ge \frac {2q^2}q=2q.
\]

Therefore, we may assume that there exists $z\in{\mathbb F}_q^2$ such that all lines through
$z$ are mapped to lines by $f$. Using suitable affine transformations, we may further assume
that 
\begin{enumerate}
\item[(i)]
$f$ maps all horizontal (vertical) lines to horizontal (vertical) lines,

\item[(ii)]
$f(0,0)=(0,0)$, $f(1,0)=(1,0)$, $f(0,1)=(0,1)$,

\item[(iii)] all lines through $(0,0)$ are mapped to lines.
\end{enumerate}

Let
\[
\begin{array}{rl}
f(x,0)=(\phi(x),0)&\quad\forall x\in{\mathbb F}_q,\cr
f(0,y)=(0,\psi(y))&\quad\forall y\in{\mathbb F}_q,\cr
\end{array}
\]
where $\phi$ and $\psi$ are permutations of ${\mathbb F}_q$ with $\phi(0)=\psi(0)=0$ and
$\phi(1)=\psi(1)=1$. For any $(x,y)\in{\mathbb F}_q^2$, it is the intersection of the vertical
line through $(x,0)$ and the horizontal line through $(0,y)$. Hence $f(x,y)$ is the
intersection  of the vertical
line through $(\phi(x),0)$ and the horizontal line through $(0,\psi(y))$, i.e.,
\[
f(x,y)=\bigl(\phi(x),\,\psi(y)\bigr).
\]
Since $f(1,1)=(1,1)$, by (iii), the line $\{(x,x):x\in{\mathbb F}_q\}$ is mapped to itself.
Hence $\phi=\psi$. Let $k\in{\mathbb F}_q$. By (iii), $f\bigl(\{(x,kx):x\in{\mathbb F}_q\}\bigr)=
\{(\phi(x),\phi(kx)):x\in{\mathbb F}_q\}$ is a line. Hence
\[
\frac{\phi(kx)}{\phi(x)}=g(k)\quad\forall x\in{\mathbb F}_q\setminus\{0\},
\]
where $g(k)\in{\mathbb F}_q$ is a function of $k$. Setting $x=1$, we have $g(k)=\phi(k)$. Thus
\begin{equation}\label{7.1}
\phi(kx)=\phi(k)\phi(x)\quad\text{for all}\ x,k\in{\mathbb F}_q.
\end{equation}

Assume that for every $a\in{\mathbb F}_q\setminus\{0\}$, all lines through $(a,0)$ which are neither horizontal nor vertical are not mapped to lines. Then we have
\[
1\text{-coaffinity}(f)\ge (q-1)(q+1-2)>2q,
\]
since $q\ge 4$.

Therefore, we may assume that there exist $a\in{\mathbb F}_q\setminus\{0\}$ and a line $L$
through $(a,0)$ such that $L$ is neither horizontal nor vertical and $f(L)$ is a line.
Let $(0,b)\in L$, where $b\in{\mathbb F}_q\setminus\{0\}$. Then $f(L)$ is the line through
$(\phi(a),0)$ and $(0,\phi(b))$. (See Figure 1.) For each $x\in{\mathbb F}_q$, the intersection
of $L$ and the vertical line through $(x,0)$ is
\[
\bigl(x,\ -\frac ba(x-a)\bigr);
\]
the intersection
of $f(L)$ and the vertical line through $(\phi(x),0)$ is
\[
\Bigl(\phi(x),\ -\frac {\phi(b)}{\phi(a)}\bigl(\phi(x)-\phi(a)\bigr)\Bigr).
\]
By (i), we have
\[
\phi\bigl(-\frac ba(x-a) \bigr)=-\frac {\phi(b)}{\phi(a)}\bigl(\phi(x)-\phi(a)\bigr).
\]
Using (\ref{7.1}), we obtain
\begin{equation}\label{7.2}
\phi(a-x)=\phi(a)-\phi(x).
\end{equation}
For any $b, x\in{\mathbb F}_q$, by (\ref{7.1}) and (\ref{7.2}),
\begin{equation}\label{7.3}
\phi(ba-bx)=\phi(b)\phi(a-x)=\phi(b)\bigl(\phi(a)-\phi(x)\bigr)=\phi(ba)-\phi(bx).
\end{equation}
Combining (\ref{7.1}) and (\ref{7.3}), $\phi$ is an automorphism of ${\mathbb F}_q$.
Hence $f\in{\rm A\Gamma L}(2,{\mathbb F}_q)$, which is a contradiction.
\end{proof}

%%%%%%%%%%% Figure 1 %%%%%%%%%%%%%%%%%%%%%%%
\vskip2mm
\setlength{\unitlength}{5mm}
\[
\begin{picture}(8,7)
\put(0,1){\vector(1,0){8}}
\put(1,0){\vector(0,1){7}}
\put(0,6){\line(1,-1){6}}
\put(3,0){\line(0,1){5}}

\put(0.8,5){\makebox(0,0)[tr]{$\scriptstyle b$}}
\put(5,0.8){\makebox(0,0)[tr]{$\scriptstyle a$}}
\put(2.8,0.8){\makebox(0,0)[tr]{$\scriptstyle x$}}
\put(6,0.2){\makebox(0,0)[bl]{$\scriptstyle L$}}
\put(3.2,3){\makebox(0,0)[bl]{$\scriptstyle (x,\,\frac ba(x-a))$}}
\put(3,3){\makebox(0,0){$\scriptstyle \bullet$}}

\end{picture}
\kern2.5cm
%%%%%%%%%%%%%%%%%
\begin{picture}(8,7)
\put(0,1){\vector(1,0){8}}
\put(1,0){\vector(0,1){7}}
\put(0,6){\line(1,-1){6}}
\put(3,0){\line(0,1){5}}

\put(0.8,5){\makebox(0,0)[tr]{$\scriptstyle \phi(b)$}}
\put(5,0.8){\makebox(0,0)[tr]{$\scriptstyle \phi(a)$}}
\put(2.8,0.8){\makebox(0,0)[tr]{$\scriptstyle \phi(x)$}}
\put(6,0.2){\makebox(0,0)[bl]{$\scriptstyle f(L)$}}
\put(3.2,3){\makebox(0,0)[bl]{$\scriptstyle (\phi(x),\,\frac {\phi(b)}{\phi(a)}(\phi(x)-\phi(a)))$}}
\put(3,3){\makebox(0,0){$\scriptstyle \bullet$}}

\end{picture}
\]

\[
\text{Figure 1. Proof of Theorem~\ref{T7.1}}
\]
\vskip2mm
%%%%%%%%%%%%%% end Figure 1 %%%%%%%%%%%%%%%%%%%%%%%%%%%%%%%%%%%%

\begin{thm}\label{T7.2}
Let $f\in {\rm Per}({\mathbb F}_3^2)\setminus{\rm AGL}(2,{\mathbb F}_3)$. Then
\[
1\text{\rm -coaffinity}(f)\ge 6.
\]
The equality holds if and only if $f\in{\rm AGL}(2,{\mathbb F}_3)\circ\tau\circ {\rm AGL}(2,{\mathbb F}_3)$, where
$\tau\in {\rm Per}({\mathbb F}_3^n)$ is any transposition.
\end{thm}

\begin{proof}
{\bf Case 1}. There do not exist two nonparallel lines in ${\mathbb F}_3^2$ which are mapped into lines by $f$. Then
\[
1\text{-coaffinity}(f)\ge 3\bigl(\frac{3^2-1}{3-1}-1\bigr)=9>6.
\]

{\bf Case 2}. There are two nonparallel lines $L_1$ and $L_2$ in ${\mathbb F}_3^2$ which are mapped into lines. 
Note that $f$ is affine
on each of these two lines. (This a special property of ${\mathbb F}_3$.) 
Therefore, through suitable affine transformations,
we may assume that $L_1={\mathbb F}_3\times\{0\}$, $L_2=\{0\}\times{\mathbb F}_3$ and
that $f|_{L_1}={\rm id}$, $f|_{L_2}={\rm id}$.

{\bf Case 2.1}. $f$ moves exactly two elements in ${\mathbb F}_3^2\setminus(L_1\cup L_2)=\{1,-1\}\times\{1,-1\}$.
$f$ is a transposition and $1\text{-coaffinity}(f)=6$.

{\bf Case 2.2}. $f$ moves exactly three elements in $\{1,-1\}\times\{1,-1\}$, say, $f=(a_1,a_2,a_3)$ where
$\{1,-1\}\times\{1,-1\}=\{a_1,a_2,a_3,a_4\}$. For each $a\in\{1,-1\}\times\{1,-1\}$, let $L_a$ be the unique line through $a$ such that $|L_a\cap(L_1\cup L_2)|=2$. Then $f(L_{a_i})$ ($1\le i\le 3$) is not a line since $a_i$ is not fixed by $f$ but the other two points on $L_{a_i}$
are. Note that the third point on the line $\overline{a_4a_i}$ ($1\le i\le 3$) is on $L_1\cup L_2$. Thus $f(\overline{a_4a_i})$ ($1\le i\le 3$) is not a line.
We also claim that for $1\le i<j\le 3$, $f(\overline{a_ia_j})$ is not a line.
Otherwise, without loss of generality, assume that $f(\overline{a_1a_2})$ is a line. $\overline{a_1a_2}$ must intersect $L_1\cup L_2$, say, at $a_0$. Note that $a_0,a_2\in f(\overline{a_1a_2})$. Hence $f(\overline{a_1a_2})
=\overline{a_0a_2}=\overline{a_0a_1}$. Thus $a_2,a_3\in f(\overline{a_1a_2})
=\overline{a_0a_1}$, which is impossible. 

Therefore,
\[
1\text{-coaffinity}(f)\ge 3+\binom 4 2 =9>6.
\]

{\bf Case 2.3}. $f=(a_1,a_2,a_3,a_4)$, where $\{1,-1\}\times\{1,-1\}=\{a_1,a_2,a_3,a_4\}$. By the argument in Cases
2.2, $f(L_{a_i})$ ($1\le i\le 4$) and $f(\overline{a_1a_2})$,  $f(\overline{a_2a_3})$, $f(\overline{a_3a_4})$,
$f(\overline{a_4a_1})$ are not lines. Hence
\[
1\text{-coaffinity}(f)\ge 8>6.
\]

{\bf Case 2.4}. $f=(a_1,a_2)(a_3,a_4)$, where $\{1,-1\}\times\{1,-1\}=\{a_1,a_2,a_3,a_4\}$. Assume $a_1=(1,1)$.
If $a_2=(1,-1)$, let $g\in{\rm GL}(2,{\mathbb F}_3)$ be given by
\[ 
g(x,y)=(x,-y).
\]
Then it is easy to see that $g\circ f$ is the transposition which moves $(0,1)$ and $(0,-1)$.

If $a_2=(-1,1)$, let $h\in{\rm GL}(2,{\mathbb F}_3)$ be given by
\[ 
h(x,y)=(-x,y).
\]
Then $h\circ f$ is the transposition which moves $(1,0)$ and $(-1,0)$.

If $a_2=(-1,-1)$, then $f(L_{a_i})$ ($1\le i\le 4$) and $f(\overline{a_1a_3})$,  $f(\overline{a_1a_4})$, $f(\overline{a_2a_3})$,
$f(\overline{a_2a_4})$ are not lines. Hence
\[
1\text{-coaffinity}(f)\ge 8>6.
\]
\end{proof}

We now turn to the proof of the Threshold Conjecture with $k=1$, $q>2$ and $n\ge 3$.

\begin{lem}\label{L7.3} 
Let $m\ge n\ge 2$ and let $f:{\mathbb F}_q^n\rightarrow {\mathbb F}_q^m$ be a one-to-one mapping which is not semi-affine. Let $1\text{\rm-coaffinity}(f)$ denote the number of lines $L$ in
${\mathbb F}_q^n$ such that $f(L)$ is not a line in ${\mathbb F}_q^m$. 
Then
\[
1\text{\rm-coaffinity}(f)\ge\frac{q^n-1}{q-1}.
\]
\end{lem}

\begin{proof}
Use induction on $n$. First assume $n=2$.

{\bf Case 1}.
There do not exist two nonparallel lines in ${\mathbb F}_q^2$ which are not mapped into line by $f$. Then
\[
1\text{-coaffinity}(f)\ge q(q+1-1)>q+1.
\]

{\bf Case 2}.
There are two nonparallel lines $L_1$ and $L_2$ in ${\mathbb F}_q^2$ such that $f(L_1)$, $f(L_2)$ are lines in
${\mathbb F}_q^m$. Since $f(L_1)$ and $f(L_2)$ are intersecting lines in ${\mathbb F}_q^m$, their affine span in
${\mathbb F}_q^m$ is a 2-flat which, without loss of generality, is assumed to be ${\mathbb F}_q^2
\times \{0\}$. 

If for every line $L_3$ in ${\mathbb F}_q^2$ such that $L_i\cap L_j$ ($1\le i<j\le 3$) are 3 distinct points $f(L_3)$ is not a line, then 
\[
1\text{-coaffinity}(f)\ge (q-1)^2\ge q+1.
\]
So, we assume that there exists a line $L_3$ in ${\mathbb F}_q^2$ such that $L_i\cap L_j$ ($1\le i<j\le 3$) are 3 distinct points and $f(L_3)$ is a line in ${\mathbb F}_q^m$.
Since $f(L_1)\cup f(L_2)\subset{\mathbb F}_q^2\times\{0\}$, it follows that
$f(L_3)\subset{\mathbb F}_q^2\times\{0\}$.

If $f({\mathbb F}_q^2)\subset {\mathbb F}_q^2\times\{0\}$, by Theorems~\ref{T7.1} and \ref{T7.2},
\[
1\text{-coaffinity}(f)\ge 2q>q+1.
\]
So we assume that $f({\mathbb F}_q^2)\not\subset{\mathbb F}_q^2\times\{0\}$.

Let $a\in {\mathbb F}_q^2\setminus(L_1\cup L_2\cup L_3)$ such that $f(a)\notin
{\mathbb F}_q^2\times\{0\}$. If $L$ is a line through $a$ such that $|L\cap(L_1\cup L_2\cup L_3)|
\ge 2$, then $f(L)$ is not a line since two points on $L$ (belonging to
$L\cap(L_1\cup L_2\cup L_3)$) are mapped to ${\mathbb F}_q^2\times\{0\}$ by $f$ but $a$ is not
mapped to ${\mathbb F}_q^2\times\{0\}$. There are only 3 lines $L_i'$ ($i=1,2,3$) in ${\mathbb F}_q^2$ such that $|L_i'\cap(L_1\cup L_2\cup L_3)|<2$: the lines through $L_j\cap L_k$ and 
parallel to $L_i$ where $\{i,j,k\}=\{1,2,3\}$.

If there are two points $a_1,a_2\in {\mathbb F}_q^2\setminus(L_1\cup L_2\cup L_3)$ such that
$f(a_i)\notin {\mathbb F}_q^2\times\{0\}$, $i=1,2$, then for any line $L$ in ${\mathbb F}_q^2$
through $a_1$ or $a_2$ with $L\ne L_i'$, $i=1,2,3$, $f(L)$ is not a line. Hence
\[
1\text{-coaffinity}(f)\ge 2(q+1)-1-3\ge q+1.
\]

If there is exactly one point $a\in{\mathbb F}_q^2\setminus(L_1\cup L_2\cup L_3)$
such that $f(a)\notin {\mathbb F}_q^2\times\{0\}$, then every line in ${\mathbb F}_q^2$
through $a$ is not a line in ${\mathbb F}_q^m$. Thus,
\[
1\text{-coaffinity}(f)\ge q+1.
\]

Now assume $n\ge 3$.

{\bf Case 1}. For any two nonparallel hyperplanes $H_1$ and $H_2$ in ${\mathbb F}_q^n$, $f$ is not semi-affine on
at least one of $H_1$ and $H_2$. Then $f$ is not semi-affine on at least 
\[
q\bigl(\frac{q^n-1}{q-1}-1\bigr)=\frac{q^2(q^{n-1}-1)}{q-1}
\]
hyperplanes in ${\mathbb F}_q^n$. By the induction hypothesis, at least $\frac{q^{n-1}-1}{q-1}$ lines on each of these
hyperplanes are not mapped into lines. Since each line in ${\mathbb F}_q^n$ lies in $\frac{q^{n-1}-1}{q-1}$ hyperplanes,
we have
\begin{equation}\label{7.4}
1\text{-coaffinity}(f)\ge\frac{\frac{q^2(q^{n-1}-1)}{q-1}\cdot \frac{q^{n-1}-1}{q-1}} {\frac{q^{n-1}-1}{q-1}}
=\frac{q^2(q^{n-1}-1)}{q-1}>\frac{q^n-1}{q-1}.
\end{equation}
 
{\bf Case 2}. There are two nonparallel hyperplanes $H_1$ and $H_2$ in ${\mathbb F}_q^n$ such that $f$ is semi-affine on both $H_1$ and $H_2$. Since $f(H_1)$ and $f(H_2)$ are $(n-1)$-flats in
${\mathbb F}_q^m$ whose intersection is an $(n-2)$-flat, their affine span
in ${\mathbb F}_q^m$ is an $n$-flat which, without loss of generality, is assumed to be
${\mathbb F}_q^n\times\{0\}$. Through a suitable semi-affine transformation, we may assume that 
\begin{equation}\label{7.5a}
f(x)=(x,0)\quad\text{for all}\ x\in H_1.
\end{equation} 
Since $f(x)=(x,0)$ for all $x\in H_1\cap H_2$ and since $\dim (H_1\cap H_2)\ge 1$, $f|_{H_2}$ must
be affine. Then it is clear that through an additional affine transformation, we may assume that
in addition to (\ref{7.5a}), 
\begin{equation}\label{7.5b}
f(x)=(x,0)\quad\text{for all}\ x\in H_2.
\end{equation} 
Thus we have
\begin{equation}\label{7.5}
f(x)=(x,0)\quad\text{for all}\ x\in H_1\cup H_2.
\end{equation} 

{\bf Case 2.1}. For every hyperplane $H_3$ in ${\mathbb F}_q^n$ such that $H_i\cap H_j$
($1\le i<j\le 3$) are 3 distinct $(n-2)$-flats, $f$ is not semi-affine on $H_3$.
By the induction hypothesis, at least $\frac{q^{n-1}-1}{q-1}$ lines on $H_3$ are not mapped into lines. The number of such hyperplanes $H_3$ is $q(\frac{q^n-1}{q-1}-2)-(q-1)$. By the same argument for (\ref{7.4}), we have
\begin{equation}\label{case2.1}
1\text{-coaffinity}(f)\ge q(\frac{q^n-1}{q-1}-2)-(q-1)>\frac{q^n-1}{q-1}.
\end{equation}

{\bf Case 2.2}. There exists a hyperplane $H_3$ in ${\mathbb F}_q^n$ such that $H_i\cap H_j$
($1\le i<j\le 3$) are 3 distinct $(n-2)$-flats and $f$ is semi-affine on $H_3$.
By (\ref{7.5}), 
\[
f(x)=(x,0)\quad\text{for all}\ x\in(H_3\cap H_1)\cup (H_3\cap H_2).
\]
Since $H_3\cap H_1$ and $H_3\cap H_2$ affinely span $H_3$, we have
\[
f(x)=(x,0)\quad\text{for all}\ x\in H_3.
\]
Thus
\[
f(x)=(x,0)\quad\text{for all}\ x\in H_1\cup H_2\cup H_3.
\]

Since $f$ is not semi-affine, $f(a)\ne (a,0)$ for some $a\in{\mathbb F}_q\setminus(H_1\cup H_2\cup H_3)$. If $L$ is a line through $a$ such that $|L\cap(H_1\cup H_2\cup H_3)|\ge 2$
and $f(a)\notin L\times\{0\}$, then $f(L)$ is not a line. (Let $b_1,b_2\in
L\cap(H_1\cup H_2\cup H_3)$. Then $f(b_i)=(b_i,0)\in L\times\{0\}$, but $f(a)\notin L\times\{0\}$.) We claim that number of lines $L$ in ${\mathbb F}_q^n$ through $a$ with
$|L\cap(H_1\cup H_2\cup H_3)|\le 1$ is
\begin{equation}\label{7.6}
7\cdot q^{n-3}+\frac{q^{n-3}-1}{q-1}.
\end{equation}
In fact, we may assume
\[
\begin{split}
H_1=\,&\{0\}\times{\mathbb F}_q\times {\mathbb F}_q\times {\mathbb F}_q^{n-3},\cr
H_2=\,&{\mathbb F}_q\times \{0\}\times{\mathbb F}_q\times {\mathbb F}_q^{n-3},\cr
H_3=\,&{\mathbb F}_q\times{\mathbb F}_q\times \{0\}\times {\mathbb F}_q^{n-3}.\cr
\end{split}
\]
Let $a=(a^{(1)},a^{(2)},a^{(3)},\dots)$ where $a^{(1)},a^{(2)},a^{(3)}\in{\mathbb F}_q\setminus\{0\}$ and let $L$ be a line through $a$ with direction vector
$v=(v^{(1)},v^{(2)},v^{(3)},\dots)$. Then $|L\cap(H_1\cup H_2\cup H_3)|\le 1$ if and only if one of the following is true:
\begin{enumerate}
\item[(i)]
at least two of $v^{(1)},v^{(2)},v^{(3)}$ are 0;

\item[(ii)]
$v^{(k)}=0$ and $(v^{(i)},v^{(j)})=t(a^{(i)},a^{(j)})$ for some $t\in{\mathbb F}_q\setminus\{0\}$
where $\{i,j,k\}=\{1,2,3\}$;

\item[(iii)]
$(v^{(1)},v^{(2)},v^{(3)})=t(a^{(1)},a^{(2)},a^{(3)})$ for some $t\in{\mathbb F}_q\setminus\{0\}$.
\end{enumerate}
Formula (\ref{7.6}) follows from these conditions.

Therefore, the number of lines $L$ in ${\mathbb F}_q^n$ through $a$ such that
$|L\cap(H_1\cup H_2\cup H_3)|\ge 2$ and $f(a)\notin L\times\{0\}$ is at least
\begin{equation}\label{7.7}
\frac{q^n-1}{q-1}-7q^{n-3}-\frac{q^{n-3}-1}{q-1}-1=q^{n-3}(q^2+q-6)-1.
\end{equation}

First assume that there are at least 4 points $a_1,\dots, a_4\in{\mathbb F}_q^n\setminus
(H_1\cup H_2\cup H_3)$ such that $f(a_i)\ne (a_i,0)$. From the above,
\[
\begin{split}
1\text{-coaffinity}(f)\,&\ge 4\bigl[q^{n-3}(q^2+q-6)-1\bigr]-\binom 42\cr
&=4q^{n-3}(q^2+q-6)-10\cr
&>\frac{q^n-1}{q-1}.\cr
\end{split}
\]

Next, assume that there are exactly $s$ points $a_1,\dots, a_s\in
{\mathbb F}_q^n\setminus(H_1\cup H_2\cup H_3)$, where $s=2$ or 3, such that $f(a_i)\ne (a_i,0)$, $1\le i\le s$. Then for every line $L$ passing through exactly one of $a_1,\dots, a_s$, $f(L)$ is not a line. Hence
\[
1\text{-coaffinity}(f)\ge s\cdot\frac{q^n-1}{q-1}-2\binom s2> \frac{q^n-1}{q-1}.
\]

Finally, assume that there is exactly one point $a\in
{\mathbb F}_q^n\setminus(H_1\cup H_2\cup H_3)$ such that $f(a)\ne(a,0)$. Then the lines in
${\mathbb F}_q^n$ which are not mapped into lines by $f$ are precisely the ones passing
through $a$. Thus,
\[
1\text{-coaffinity}(f)=\frac{q^n-1}{q-1}.
\]
\end{proof}
 
\begin{thm}\label{T7.4}
Let $n\ge 3$ and 
$f\in {\rm Per}({\mathbb F}_q^n)\setminus {\rm A\Gamma L}(n,{\mathbb F}_q)$. Then
\[
1\text{\rm -coaffinity}(f)\ge 2q\genfrac{[}{]}{0pt}{0}{n-1} 1_q=\frac{2q(q^{n-1}-1)}{q-1}.
\]
The equality holds if and only if $f\in {\rm A\Gamma L}(n,{\mathbb F}_q)
\circ\tau\circ {\rm A\Gamma L}(n,{\mathbb F}_q)$, where $\tau\in {\rm Per}({\mathbb F}_q^n)$ is any transposition.
\end{thm}

\begin{proof}
The arguments in this proof are very similar to those in the proof of Lemma~\ref{L7.3}.

{\bf Case 1}.
For any two nonparallel hyperplanes $H_1$ and $H_2$ in ${\mathbb F}_q^n$, $f$
is not semi-affine on at least one of $H_1$ and $H_2$. By Lemma~\ref{L7.3}
and (\ref{7.4}),
\[
1\text{-coaffinity}(f)\ge\frac{q^2(q^{n-1}-1)}{q-1}>\frac{2q(q^{n-1}-1)}{q-1}.
\]

{\bf Case 2}. 
There are two nonparallel hyperplanes $H_1$ and $H_2$ in ${\mathbb F}_q^n$ such that $f$ is semi-affine on both $H_1$ and $H_2$. By the proof of
Lemma~\ref{L7.3}, we may assume that  $f|_{H_1}=$ id, $f|_{H_2}=$ id.

{\bf Case 2.1}.
For every hyperplane $H_3$ in ${\mathbb F}_q^n$ 
such that $H_i\cap H_j$ ($1\le i<j\le 3$) are 3 distinct $(n-2)$-flats,
$f$ is not semi-affine on $H_3$. 
By (\ref{case2.1}), we have
\[
1\text{-coaffinity}(f)\ge q\Bigl(\frac{q^n-1}{q-1}-2\Bigr)-(q-1)>\frac{2q(q^{n-1}-1)}{q-1}.
\]

{\bf Case 2.2}. There exists a hyperplane $H_3$ in ${\mathbb F}_q^n$
such that $H_i\cap H_j$ ($1\le i<j\le 3$) are 3 distinct $(n-2)$-flats and
$f$ is semi-affine on $H_3$.
By the same argument in Case 2.2 of the proof of Lemma~\ref{L7.3}, we have 
\[
f(x)=x\quad\text{for all}\ x\in H_1\cup H_2\cup H_3.
\]

First assume that there are at least 5 elements $a_1,\dots,a_5\in
{\mathbb F}_q^n\setminus(H_1\cup H_2\cup H_3)$ such that $f(a_i)\ne a_i$, $1\le i\le 5$.
By (\ref{7.7}), we have
\begin{equation}\label{7.8}
1\text{-coaffinity}(f)\ge 5\bigl[q^{n-3}(q^2+q-6)-1\bigr]-\binom 5 2=
5q^{n-3}(q^2+q-6)-15.
\end{equation}
When $q\ge 4$, we have
\[
5q^{n-3}(q^2+q-6)-15>\frac{2q(q^{n-1}-1)}{q-1}.
\]
When $q=3$, any line $L$ in ${\mathbb F}_3^n$ with $|L\cap(H_1\cup H_2\cup H_3)|\ge 2$ cannot
pass through more than one of $a_1,\dots,a_5$ since $|L|=3$. Thus (\ref{7.8}) can be
improved to
\[
1\text{-coaffinity}(f)\ge 5\bigl[q^{n-3}(q^2+q-6)-1\bigr]
>\frac{2q(q^{n-1}-1)}{q-1}.
\]

Next, assume that there are exactly $s$ elements $a_1,\dots,a_s\in
{\mathbb F}_q^n\setminus(H_1\cup H_2\cup H_3)$, where $s=3$ or 4, such that $f(a_i)\ne a_i$,
$1\le i\le s$. Then for every line $L$ passing through exactly one of
$a_1,\dots,a_s$, $f(L)$ is not a line. Hence
\[
1\text{-coaffinity}(f)\ge s\,\frac{q^n-1}{q-1}-2\binom s 2=\frac s{q-1}
\bigl[q^n-(s-1)q+s-2\bigr]
>\frac{2q(q^{n-1}-1)}{q-1}.
\]

Finally, assume that there are exactly 2 elements $a_1,a_2\in{\mathbb F}_q^n
\setminus(H_1\cup H_2\cup H_3)$ such that $f(a_i)\ne a_i$, $i=1,2$. Then $f$ is a transposition.
\end{proof} 

Combining Theorems~\ref{T7.1}, \ref{T7.2} and \ref{T7.4}, we have the following corollary.

\begin{cor}\label{C7.5}
The Threshold Conjecture holds for $q\ge 3$, $k=1$ and $n>1$.
\end{cor}

\noindent {\bf Remark}. Theorems~\ref{T7.2} and \ref{T7.4} state that for $q=3$, $k=1$, $n\ge 2$ 
or $q>3$, $k=1$, $n\ge 3$, $1\text{-coaffinity}(f)=2q\genfrac{[}{]}{0pt}{1}{n-1}1_q$ only
when $f\in{\rm A\Gamma L}(n,{\mathbb F}_q)\circ\tau\circ {\rm A\Gamma L}(n,{\mathbb F}_q)$ for some transposition $\tau\in{\rm Per}({\mathbb F}_q^n)$. It is an open question whether the same holds
for $q>3$, $k=1$ and $n=2$.

%%%%%%%%%%%%%%%%%%%%%%%%%%%%%%%%%%%%%%%%%%%%%%%%%
\section*{Acknowledgment}

The proof of Theorem~\ref{T7.1} is based on a method by Professor W-X. Ma.

\end{document}